\input amstex
\documentstyle{amsppt}
\input FIG.sty

\topmatter
\title Improving Rogers' upper bound \\ for the density of unit ball
packings \\ via estimating the surface area of Voronoi cells from below \\
in Euclidean $d-$space for all $d\ge 8$
\endtitle
\author K\'aroly Bezdek
\endauthor
\leftheadtext {K\'aroly Bezdek}
\rightheadtext{ Improving Rogers' upper bound}
\affil Department of Geometry, E\"otv\"os University, Budapest \\  and \\ 
       Department of Mathematics, Cornell University, Ithaca NY
\endaffil

\address   E\"otv\"os University, Department of Geometry,  
           H-1117 Budapest, P\'azm\'any P\'eter s\'et\'any 1/C, Hungary,  
           e-mail: kbezdek\@ludens.elte.hu and \newline
           Cornell University, Department of Mathematics, 
           310 Malott Hall, Ithaca, NY 14853-4201, USA, 
           e-mail: bezdek\@math.cornell.edu 
\endaddress

\date October 17, 2001
\enddate

\thanks  The author was partially supported by the Hung. Nat. Sci. Found.
(OTKA), grant no. T029786, and by the Combinatorial Geometry Project
of the Research Found. FKFP0151/1999.
\endthanks

\keywords  Sphere packings, Voronoi cells, surface area, density
\endkeywords

\subjclass Primary 52C17, 52A40; Secondary 52B60, 52A38
\endsubjclass

\abstract { The sphere packing problem asks for the densest packing of unit balls in ${\bold E^d}$. This
problem has its roots in geometry, number theory and information theory and it is part of Hilbert's 18th problem. One of the most attractive results on the sphere packing problem was proved by C. A. Rogers
in 1958. It can be phrased as follows. Take a regular $d-$dimensional simplex of edge length $2$ in ${\bold E^d}$ and then draw a $d-$dimensional unit ball around each vertex of the simplex. Let $\sigma _d$ denote the ratio of the volume of the portion of the simplex covered by balls to the volume of the simplex. Then the volume of any Voronoi cell in a packing of unit balls in ${\bold E^d}$ is at least $\frac{\omega_d}{\sigma_d}$, where $\omega_d$ denotes the volume of a $d-$dimensional unit ball. This has the immediate corollary that the density of any unit ball packing in ${\bold E^d}$ is at most $\sigma_d$.
In 1978 Kabatjanskii and Leven\v stein improved this bound
for large $d$. In fact, Rogers' bound is the presently known best bound for $4\le d\le 42$, and above that the Kabatjanskii-Leven\v stein bound takes over.  In this paper we improve Rogers' upper bound for the density 
of unit ball packings in Euclidean $d-$space for all $d\ge 8$ and improve Kabatjanskii-Leven\v stein upper bound in small dimensions. Namely, we show that the volume of any Voronoi cell in a packing of unit balls in ${\bold E^d}, d\ge 8$ is at least $\frac{\omega_d}{\widehat\sigma_d}$ and so the
density of any unit ball packing in ${\bold E^d}, d\ge 8$ is at most $\widehat\sigma_d$, where
$\widehat\sigma_d$ is a geometrically well-defined quantity satisfying the inequality 
$ \widehat\sigma_d<\sigma_d$ for all $d\ge 8$. We prove this by showing that the surface area of any Voronoi cell in a 
packing of unit balls in ${\bold E^d}, d\ge 8$ is at least $\frac{d\cdot\omega_d}{\widehat\sigma_d}$.}  
\endabstract
\define\deriv#1{\dfrac{d#1}{dx}}
\define\derive#1{\dfrac{d#1}{dh}}

\endtopmatter

\document

\heading 0. Introduction
\endheading

A family of non-overlapping $d-$dimensional balls of radii $1$ in the $d-$dimensional Euclidean space ${\bold E^d}$ is called a unit
ball packing of ${\bold E^d}$. The density of the packing is the proportion of space covered by these unit balls. The sphere packing
problem asks for the densest packing of unit balls in ${\bold E^d}$. Indubitably, of all problems concerning packing it was the sphere
packing problem which attracted the most attention in the past decade. It has its roots in geometry, number theory and information
theory and it is part of Hilbert's 18th problem. The reader is referred to \cite {10} (especially the third edition, which has about 800 references covering 1988-1998) for further information, definitions and references. In what follows we report on a few selected developments and then we state the main results of this paper.

The Voronoi cell of a unit ball in a packing of unit balls in ${\bold E^d}$ is the set of points that are not farther away from the
center of the given ball than from any other ball's center. As it is well-known (see for exapmle \cite{24}) the Voronoi cells of a unit
ball packing in ${\bold E^d}$ form a tiling of ${\bold E^d}$. One of the most attractive results on the sphere packing problem was
proved by C. A. Rogers \cite{23} in 1958. It was rediscovered by Baranovskii \cite {2} and extended to spherical and hyperbolic spaces by B\"or\"oczky \cite {6}. It can be phrased as follows. Take a regular $d-$dimensional simplex of edge length $2$ in ${\bold E^d}$ and then draw a $d-$dimensional unit ball around each vertex of the simplex. Let $\sigma _d$ denote the ratio of the volume of the portion of the simplex covered by balls to the volume of the simplex. Then the volume of any Voronoi cell in a packing of unit balls in ${\bold E^d}$ is at least $\frac{\omega_d}{\sigma_d}$, where $\omega_d$ denotes the volume of a $d-$dimensional unit ball. This has the immediate corollary that the (upper) density of any unit ball packing in ${\bold E^d}$ is at most $\sigma_d$.
Daniel's asymptotic formula \cite{24} yields that  
$$\sigma_d=\frac{d}{e}2^{-(0.5+o(1))d}\ (\text{\rm as \ }d\to\infty).$$
Then 20 years later, in 1978 Kabatjanskii and Leven\v stein \cite{18} improved this bound in the exponential order of magnitude as follows. They showed that the density of any unit ball packing in ${\bold E^d}$ is at most
$$2^{-(0.599+o(1))d}\ (\text{\rm as \ }d\to\infty).$$
In fact, Rogers' bound is the presently known best bound for $4\le d\le 42$ (see also Remark 2 below), and above that the
Kabatjanskii-Leven\v stein bound takes over (\cite{10}, p. 20). 

There has been some important recent progress concerning the existence of economical packings. On the one hand, improving earlier results, Ball \cite{1} proved through a very elegant completely new variational argument that for each $d$,
there is a lattice packing of unit balls in ${\bold E^d}$ with density at least
$$\frac{d-1}{2^{d-1}}\zeta (d),$$
where $\zeta (d)=\sum_{k=1}^{\infty}\frac{1}{k^d}$ is the Riemann zeta function. On the other hand, 
for some small values of $d$, there are 
explicit (lattice) packings which give densities (considerably) higher than the bound just stated. In connection with this we briefly mention some of the exciting results of lower dimensions. The reader is
referred to \cite{10}, \cite{12}, \cite{22} and \cite{25} for a comprehensive view of results of this type.

In ${\bold E^3}$, the face centered cubic lattice packing of congruent balls with density $\frac{\pi}{\sqrt{18}}=0.74048\dots$ was
conjectured by Kepler to be the densest packing among all packings of unit balls  (\cite{10}). It has been proven that the face
centered cubic lattice is a locally optimal arrangement (for completely independent approaches to this see \cite{3}, \cite{11} and
\cite{13}) and that the density of any packing of congruent balls in ${\bold E^3}$ cannot be larger than 0.773055\dots , which was
established by Muder in \cite{21}. Finally, in 1998 Hales \cite{16} announced that the final step in the proof of Kepler Conjecture has
been completed : the Kepler conjecture is now a theorem (finishing a several years complex project that started with \cite{14} and
\cite{15}). Hales' proof is long and difficult and requires extensive computer calculation. As of October, 2001, it has not yet been
published but, it is widely regarded as being likely to be correct. In ${\bold E^d, 4\le d\le 9}$,
the densest packings of congruent balls known are obtained by the "laminating" or "greedy" construction
described in \cite{8}. (In fact, \cite{8} determines all equivalent laminated lattices for $d\le 25$ and
produces many of the densest lattices known up to 25.) In ${\bold E^{10}}$, we encounter for the first time
a nonlattice packing of congruent balls that is denser than all known lattice packings. This packing is obtained
from "Construction $A$" (for more details see \cite{20} and \cite{25}). In ${\bold E^d}, 18\le d\le 22$, record nonlattice packings of
congruent balls have recently been given such as Vardy's construction \cite{26} ("Construction $B^*$") and Bierbrauer and Edel record
packing in ${\bold E^{18}}$ (see \cite{25}
and also \cite{9}). The packing of congruent balls with proper radii around the points of the Leech lattice is a remarkably dense packing in ${\bold E^{24}}$. New packings of congruent balls in ${\bold E^d},
26\le d\le 31$ have been recently discovered by Bacher, Borcherds, Conway, Sloane, Vardy, Venkov - see \cite{10}
for details. Finally, we mention the Kschischang and Pasupathy lattice packing \cite{19} of balls in ${\bold E^{36}}$ and the Mordell-Weil lattice packings of balls in ${\bold E^{d}}, 80\le d\le 4096$ discovered by
Elkies and Shioda that are the densest lattice packings of balls known in those dimensions -see 
\cite{10} and \cite{25} for details. For a complete account on record packings we refer to \cite{22}. All these explicit constructions raise the well-known challenging question whether one can find a smaller upper bound than Rogers' bound for the density of unit ball packings, especially in low dimensions. The Corollary established in this paper does exactly this by improving Rogers' upper bound for the density of unit ball packings in Euclidean $d-$space for all $d\ge 8$.

In what follows we state the major results of this paper. As usual, let $\text{\rm lin}(\dots ),$
$\text{\rm aff}(\dots ),$ $\text{\rm conv}(\dots ),$ $\text{\rm Vol}_d(\dots ),$ $\omega_d,$ 
$\text{\rm SVol}_{d-1}(\dots ),$ $\text{\rm dist}(\dots ),$ $\Vert\dots \Vert$ and ${\bold o}$ refer to the linear hull, the affine hull, the convex hull in ${\bold E^d}$, the $d-$dimensional Euclidean volume measure, the $d-$dimensional volume of a $d-$dimensional unit ball, the $(d-1)-$dimensional spherical volume measure, the distance function in ${\bold E^d}$, the standard Euclidean norm and to the origin in ${\bold E^d}$.

Let $\text{\rm conv}\{{\bold o}, {\bold w}_1, \dots , {\bold w}_d\}$ be a $d-$dimensional simplex having the property that the linear hull $\text{\rm lin}\{{\bold w}_j-{\bold w}_i\vert i< j\le d\}$ is orthogonal to
the vector ${\bold w}_i$ in ${\bold E^d}, d\ge 8$ for all $1\le i\le d-1$ that is let 
$$\text{\rm conv}\{{\bold o}, {\bold w}_1, \dots , {\bold w}_d\}$$ 
be a $d-$dimensional orthoscheme in ${\bold E^d}$ moreover, let  
$$\Vert {\bold w}_i\Vert=\sqrt{\frac{2i}{i+1}} \text{\rm \ for \ all \ }1\le i\le d.$$
It is clear that in the right triangle $\vartriangle {\bold w}_{d-2}{\bold w}_{d-1}{\bold w}_{d}$ with
right angle at the vertex ${\bold w}_{d-1}$ we have the inequality $\Vert {\bold w}_{d}-{\bold w}_{d-1}\Vert=
\sqrt{\frac{2}{d(d+1)}}<\sqrt{\frac{2}{(d-1)d}}=\Vert{\bold w}_{d-1}-{\bold w}_{d-2}\Vert$ and therefore
$\angle {\bold w}_{d-1}{\bold w}_{d-2}{\bold w}_{d}<\frac{\pi}{4}$ (see Fig. 0). Now, in the plane 
$\text{\rm aff}\{{\bold w}_{d-2}, {\bold w}_{d-1},$ ${\bold w}_{d}\}$ of the triangle
$\vartriangle {\bold w}_{d-2}{\bold w}_{d-1}{\bold w}_{d}$ let 
$$\vartriangleleft {\bold w}_{d-2}{\bold w}_{d}{\bold w}_{d+1}$$ 
denote the circular sector of central angle $\angle {\bold w}_{d}{\bold w}_{d-2}{\bold w}_{d+1}=
\frac{\pi}{4}-\angle {\bold w}_{d-1}{\bold w}_{d-2}{\bold w}_{d}$ and of center ${\bold w}_{d-2}$ sitting over
the circular arc with endpoints ${\bold w}_{d}, {\bold w}_{d+1}$ and radius $\Vert {\bold w}_d-{\bold w}_{d-2}
\Vert= \Vert{\bold w}_{d+1}-{\bold w}_{d-2}\Vert$ such that $\vartriangleleft {\bold w}_{d-2}{\bold w}_{d}
{\bold w}_{d+1}$ and $\vartriangle {\bold w}_{d-2}{\bold w}_{d-1}{\bold w}_{d}$ are adjacent along the line
segment ${\bold w}_{d-2}{\bold w}_{d}$ and are separated by the line of ${\bold w}_{d-2}{\bold w}_{d}$.
Then let
$$D({\bold w}_{d-2},{\bold w}_{d-1},{\bold w}_{d},{\bold w}_{d+1})= \vartriangle {\bold w}_{d-2}{\bold w}_{d-1}{\bold w}_{d}\cup \vartriangleleft {\bold w}_{d-2}{\bold w}_{d}{\bold w}_{d+1}$$
be the convex domain generated by the triangle $\vartriangle {\bold w}_{d-2}{\bold w}_{d-1}{\bold w}_{d}$
with constant angle $\angle{\bold w}_{d-1}{\bold w}_{d-2}{\bold w}_{d+1}=\frac{\pi}{4}$.

\OneFigure{file=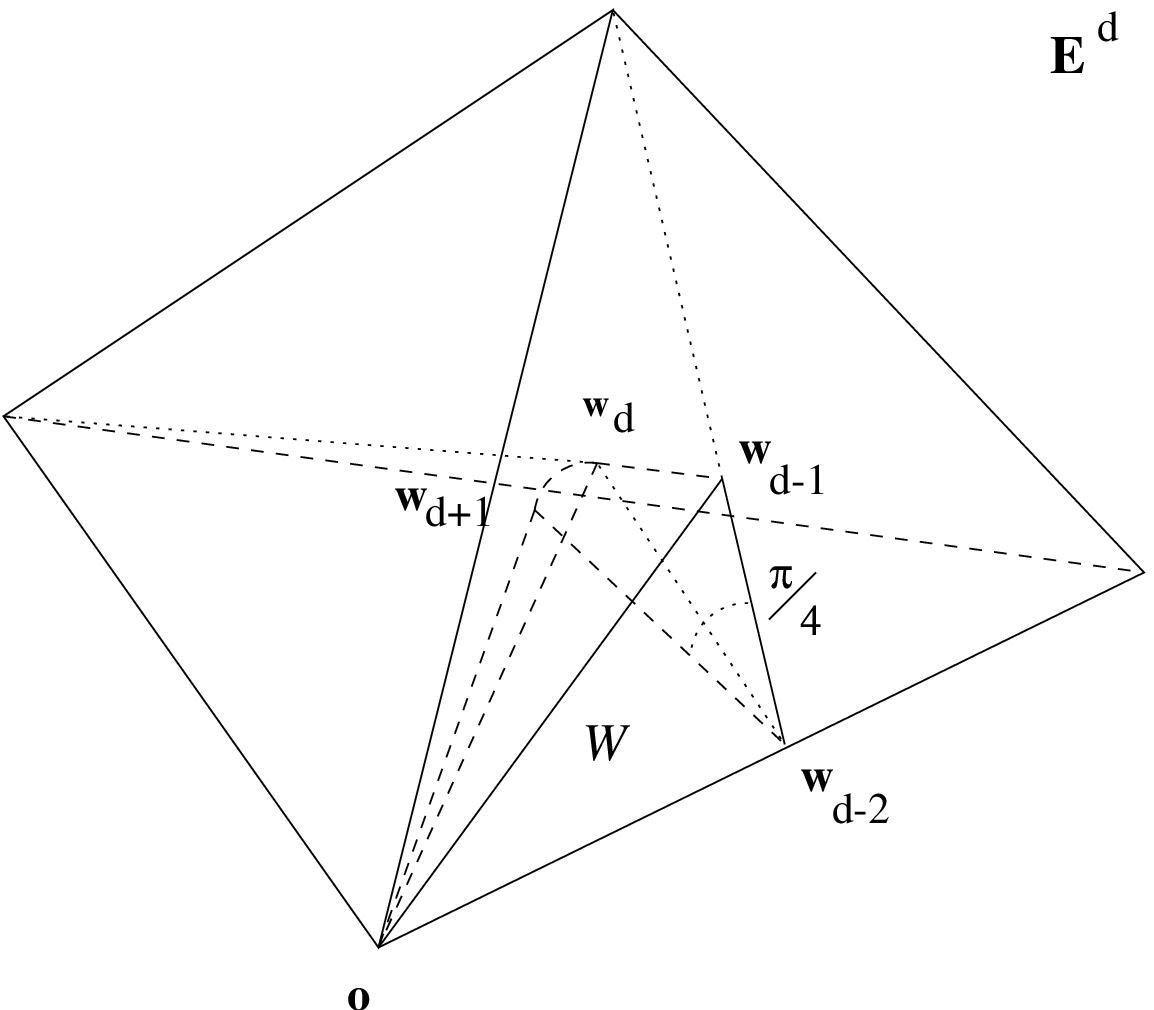}{Figure 0}

Now, let
$$W=\text{\rm conv}\left(\{{\bold o}, {\bold w}_1, \dots , {\bold w}_{d-3}\}\cup D({\bold w}_{d-2},{\bold w}_{d-1},{\bold w}_{d},{\bold w}_{d+1})\right)$$
be the $d-$dimensional wedge (or cone) with $(d-1)-$dimensional base 
$$Q_{W}=\text{\rm conv}\left(\{{\bold w}_1, \dots , {\bold w}_{d-3}\}\cup D({\bold w}_{d-2},{\bold w}_{d-1},{\bold w}_{d},{\bold w}_{d+1})\right) \text{ \ and \ apex\ }{\bold o}.$$
Finally, if $B=\{ {\bold x}\in {\bold E^d}\vert \ \text{dist}({\bold o}, {\bold x})=\Vert{\bold x}\Vert\le 1\}$ 
denotes the $d-$dimensional unit ball centered at the origin ${\bold o}$ of $\bold E^d$ and $S=\{ {\bold x}\in {\bold E^d}\vert \ \text{dist}({\bold o}, {\bold x})=\Vert{\bold x}\Vert = 1\} $ denotes the $(d-1)-$dimensional unit sphere centered at ${\bold o}$, then let 
$$\widehat{\sigma}_d= \frac{\text{\rm SVol}_{d-1}(W\cap S)}{\text{\rm Vol}_{d-1}(Q_{W})}=\frac{\text{\rm Vol}_{d}(W\cap B)}{\text{\rm Vol}_{d}(W)}$$
be the the surface density (resp., volume density) of the unit sphere $S$ (resp., of the unit ball $B$) in the wedge
$W$. For the sake of completeness we remark that as the regular $d-$dimensional simplex of edge length $2$ can be dissected into $(d+1)!$ pieces each being congruent to $ \text{\rm conv}\{{\bold o}, {\bold w}_1, \dots , {\bold w}_d\}$ therefore $$\sigma_d=\frac{\text{\rm Vol}_{d}( \text{\rm conv}\{{\bold o}, {\bold w}_1, \dots , {\bold w}_d\}  \cap B)}{\text{\rm Vol}_{d}( \text{\rm conv}\{{\bold o}, {\bold w}_1, \dots , {\bold w}_d\} )}.$$

Now, we are ready to state the main result of this paper. Recall that the surface density of any unit
sphere in its Voronoi cell in a unit sphere packing of ${\bold E^d}$ is defined as the ratio of the surface area of the unit sphere to the surface area of its Voronoi cell. The following theorem improves Theorem 3
of \cite{4} which claims that the surface area of any Voronoi cell in a packing of
unit balls in ${\bold E^d}, d\ge 2$ is at least $\frac{d\cdot \omega_d}{\sigma_d}$.

\proclaim {Theorem} The surface area of any Voronoi cell in a packing
of unit balls in the $d-$dimensional Euclidean space ${\bold E^d}, d\ge 8$ is at least $\frac{d\cdot \omega_d}{\widehat\sigma_d}$, that is the surface density of any unit sphere in its Voronoi cell in a unit sphere packing of ${\bold E^d}, d\ge 8$ is at most $\widehat\sigma_d$.
\endproclaim

As the volume of a Voronoi cell in a unit ball packing of ${\bold E^d}$ is at least as large as
$\frac{1}{d}$ times the surface area of the Voronoi cell the following result follows from the Theorem.

\proclaim {Corollary} The volume of any Voronoi cell in a packing of
unit balls in ${\bold E^d}, d\ge 8$ is at least $\frac{\omega_d}{\widehat\sigma_d}$. Thus,
the (upper) density of any unit ball packing in ${\bold E^d}, d\ge 8$ is
at most $\widehat\sigma_d$.
\endproclaim

Finally, we show that our upper bound $\widehat\sigma_d $ for the density of unit ball
packings in Euclidean $d-$space is indeed better than Rogers' bound $\sigma_d$.   

\proclaim{ Proposition } $\widehat\sigma_d<\sigma_d$ for all $d\ge 8$.
\endproclaim

{\bf Remark 1.} It is not hard to see that the proof of the Theorem presented in the sections below can be used
to prove the following stronger statement. Take a Voronoi cell of a unit ball in a packing of unit balls in the
$d-$dimensional Euclidean space ${\bold E^d}, d\ge 8$ and then take the intersection of the given Voronoi cell with
the closed $d-$dimensional ball of radius $\sqrt{\frac{2d}{d+1}}$ concentric to the unit ball of the Voronoi cell. Then the
surface area of the truncated Voronoi cell is at least $\frac{d\cdot \omega_d}{\widehat\sigma_d}$.
\bigskip
{\bf Remark 2.} Our density bound $\widehat\sigma_d$ for the density of unit ball packings in ${\bold E^d}, d\ge 8$
is a geometrically explicit bound for all $d\ge 8$. However, in concrete small dimensions in particular, in dimension $8$
our method of proof suggests further improvements via the $3-$dimensional skeleton of Voronoi polytopes. As a
result the best possible numerical value for the density bound generated by our method in dimension $8$ is still in progress.
In connection with this we mention that Cohn and Elkies just recently announced an analogue for sphere packing of the linear
programming bounds for error-correcting codes in \cite{7} and used it to improve Rogers' upper bound for the density of sphere packings
for dimensions $4$ through $36$. This work of Cohn and Elkies seems to be in progress as well. Finally, W.-Y. Hsiang \cite{17} 
very recently announced a solution of the $8-$dimensional sphere packing problem but, details are not yet public. Their methods 
are appearantly quite different from ours.
\bigskip

The organization of the rest of the paper is as follows. The proof of the Theorem consists of 6 major steps
which are discussed in consecutive separate sections of the paper. Then a short proof of the 
Proposition is presented in the last section. Finally, we have to emphasize that although most of the 
methods of this paper work in all dimensions being at least $8$ they are designed to act together in
an efficient way mostly in low dimensions.

\heading 1. Some metric properties of Voronoi cells \\  
         of unit ball packings in ${\bold E^d}$
\endheading

Let $P$ be a bounded Voronoi cell i.e. a $d-$dimensional Voronoi polytope of
a packing $\Cal P$ of $d-$dimensional unit balls in $\bold E^d$. Without loss of generality we may assume that the unit ball $B=\{ {\bold x}\in {\bold E^d}\vert \ \text{dist}({\bold o}, {\bold x})=\Vert{\bold x}\Vert\le 1\}$ centered at the origin ${\bold o}$ of $\bold E^d$ is one of the unit balls of $\Cal P$ with $P$ as its Voronoi cell. Then $P$ is the intersection of
finitely many closed halfspaces of $\bold E^d$ each of which is bounded
by a hyperplane which is the perpendicular bisector of a line segment
${\bold o}{\bold x}$ with ${\bold x}$ being the center of some unit ball
of $\Cal P$. Now, let $F_{d-i}$ be an arbitrary $(d-i)-$dimensional face of $P$, $1\le i\le d$. Then clearly there are at least 
$i+1$ Voronoi cells of $\Cal P$ which meet along the face $F_{d-i}$ 
that is contain $F_{d-i}$ (one of which is of course, $P$). Also, it is clear from the construction that the
affine hull of centers of the unit balls sitting in all of these Voronoi cells is orthogonal to $\text{aff}F_{d-i}$. Thus, there are unit balls of these
Voronoi cells with centers $\{{\bold o}, {\bold x}_1, \dots ,{\bold x}_i\}$ 
such that $X=\text{conv}\{{\bold o}, {\bold x}_1, \dots ,
{\bold x}_i\}$ is an $i-$dimensional simplex and of course, $\text{aff}X$ is orthogonal to $\text{aff} F_{d-i}$. Hence, if $R(F_{d-i})$ denotes the radius of the $(i-1)-$dimensional sphere that passes through the vertices of $X$, then $$R(F_{d-i})=\text{\rm dist}({\bold o}, \text{\rm aff} F_{d-i}),
\text{\rm \ where\ }1\le i\le d.$$

\proclaim {Lemma 1} If $F_{d-i-1}\subset F_{d-i}$ and $R(F_{d-i})=R < \sqrt{2}$ for some $i, 1\le i\le d-1$,
then $$ \frac{2}{\sqrt{4-R^2}}\le R(F_{d-i-1}).$$
\endproclaim
\subhead Proof
\endsubhead
Recall that if $C$ is the convex hull of the centers of the unit balls 
of ${\Cal P}$ whose Voronoi cells contain $F_{d-i}$,
then $\text{\rm aff}C$ is an $i-$dimensional affine subspace being totally
orthogonal to $\text{\rm aff}F_{d-i}$ moreover, $\text{\rm aff}F_{d-i}\cap C={\bold b}$ is the center of the $(i-1)-$dimensional sphere that passes
through the vertices of $C$. As ${\bold o}$ is among the vertices of $C$ there
exist vertices of $C$ say, ${\bold c}_1, \dots , {\bold c}_i$ such that
$C_i=\text{\rm conv}\{{\bold o}, {\bold c}_1, \dots , {\bold c}_i\}$ is an
$i-$dimensional simplex containing ${\bold b}$ with
$$\Vert {\bold b}-{\bold o}\Vert =\Vert {\bold b}-{\bold c}_1\Vert =
\dots =\Vert {\bold b}-{\bold c}_i\Vert =R(F_{d-i})=R.$$ 

From now on we deal
only with the case 
$${\bold b}\in \text{\rm relint}C_i$$
for the reason that the case ${\bold b}\in\text{\rm relbd}C_i$ follows from this by standard limit procedure. 

\OneFigure{file=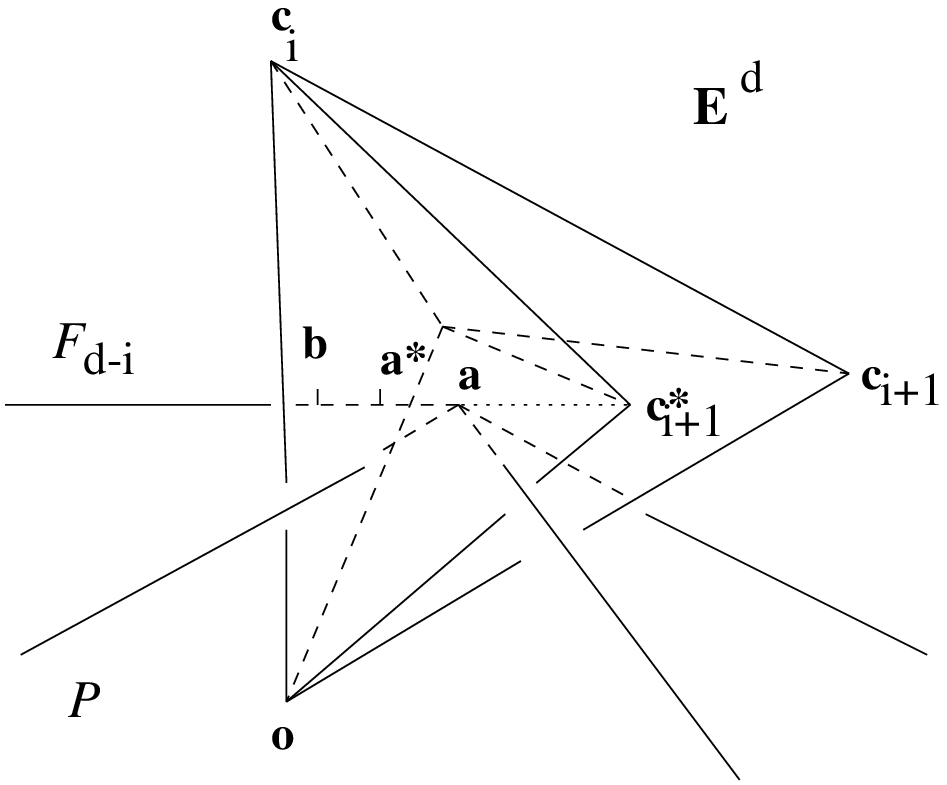}{Figure 1}

As $F_{d-i-1}\subset F_{d-i}$ there must be a unit ball of $\Cal P$ with center
say, ${\bold c}_{i+1}$ whose Voronoi cell intersects $F_{d-i}$ precisely in $F_{d-i-1}$. Thus, $$C_{i+1}=\text{\rm conv}\{{\bold o}, {\bold c}_1, \dots ,
{\bold c}_i, {\bold c}_{i+1}\}$$ is an $(i+1)-$dimensional simplex moreover,
if the center of the $i-$dimensional sphere that passes through the vertices
of $C_{i+1}$  is ${\bold a}$, then 
$$\Vert {\bold a}-{\bold o}\Vert = \Vert {\bold a}-{\bold c}_1\Vert = \dots =
\Vert {\bold a}- {\bold c}_i\Vert =\Vert {\bold a}- {\bold c}_{i+1}\Vert =
R(F_{d-i-1}).$$
Obviously, all the edges of $C_{i+1}$ have lengths larger than or equal to $2$
(see Fig. 1).


Finally, let ${\bold c}_{i+1}^*\in {\bold E^d}$ be a point such that
$$\Vert {\bold c}_{i+1}^*-{\bold o}\Vert =\Vert {\bold c}_{i+1}^*-{\bold c}_1
\Vert =\dots =\Vert {\bold c}_{i+1}^*-{\bold c}_i\Vert =2.$$
Then $$C_{i+1}^*=\text{\rm conv}\{{\bold o}, {\bold c}_1, \dots , {\bold c}_i,
{\bold c}_{i+1}^*\}$$ is an $(i+1)-$dimensional simplex moreover, if ${\bold a}^*$ is the center of the $i-$dimensional sphere that passes through the
vertices of $C_{i+1}^*$, then 
$${\bold a}^*\in \text{\rm relint}C_{i+1}^*\text{\rm \ and\ }$$
$$\Vert{\bold a}^*-{\bold o}\Vert =\Vert{\bold a}^*-{\bold c}_1\Vert =
\dots =\Vert {\bold a}^*- {\bold c}_i\Vert =\Vert {\bold a}^* -
{\bold c}_{i+1}^*\Vert =\frac{2}{\sqrt{4-R^2}}.$$
Notice that the length of each edge of $C_{i+1}$ is larger than or equal
to the length of the corresponding edge of $C_{i+1}^*$.

Now, recall the following statement which is a special case of 
Lemma 7 in \cite{5}.
\proclaim {Sublemma 1} Let $Y^*\subset {\bold E^n}, 1\le n$ be an $n-$dimensional simplex with vertices $\{{\bold y}_0^*, {\bold y}_1^*, \dots , {\bold y}_n^*\}\subset
{\bold E^n}$ and let ${\bold p}^*\in \text{\rm int}Y^*$ be an arbitrary 
interior point of $Y^*$.
If $Y\subset {\bold E^n}$ is an $n-$dimensional simplex with vertices
$\{{\bold y}_0, {\bold y}_1, \dots , {\bold y}_n\}\subset {\bold E^n}$ 
and ${\bold p}\in {\bold E^n}$ is a point with the property that
$\Vert{\bold y}_k^*-{\bold y}_l^*\Vert\le\Vert{\bold y}_k-{\bold y}_l\Vert$
for all $0\le k<l\le n$ and $\Vert{\bold p}^*-{\bold y}_j^*\Vert\ge
\Vert{\bold p}-{\bold y}_j\Vert$ for all $0\le j\le n$, then 
the simplices $Y^*$ and $Y$ are congruent moreover, $\Vert{\bold p}^*-{\bold y}_j^*\Vert = \Vert{\bold p}-{\bold y}_j\Vert$ for all $0\le j\le n$.
\endproclaim

Thus, if one assumes that $R(F_{d-i-1})<\frac{2}{\sqrt{4-R^2}}$, then
Sublemma 1 applied to the simplices $C_{i+1}^*$ and $C_{i+1}$ immediately
leads to a contradiction. This completes the proof of Lemma 1.$\hskip8.0cm \qed$

As an easy corollary of Lemma 1 we get the following well-known inequality (\cite{23}).

\proclaim {Corollary 1} $\sqrt{\frac{2i}{i+1}}\le R(F_{d-i})$ for all $1\le i\le d$.
\endproclaim
We will need the following metric property of Voronoi polytopes as well. (For
a somewhat weaker version of this see \cite{6, pages 257-258}.)

\proclaim {Lemma 2} If $R(F_{d-i})< \sqrt{2}$ for some $i, 1\le i\le d$, then
the orthogonal projection of ${\bold o}$ onto $\text{\rm aff}F_{d-i}$
belongs to $\text{\rm relint}F_{d-i}$ and so
$R(F_{d-i})=\text{\rm dist}({\bold o}, F_{d-i})$.
\endproclaim
\subhead Proof
\endsubhead
Take the $i-$dimensional simplex $C_i$ defined in the beginning of the proof of Lemma 1. The vertices of $C_i$ are ${\bold o}, {\bold c}_1, \dots , {\bold c}_i$
that are centers of unit balls of $\Cal P$ whose Voronoi cells all contain
$F_{d-i}$. Thus, $\text{\rm aff}C_i$ is orthogonal to $\text{\rm aff}F_{d-i}$
moreover,
$${\bold b}=C_i\cap\text{\rm aff}F_{d-i}\text{\rm \ and \ }
\Vert{\bold b}-{\bold o}\Vert =\Vert {\bold b}-{\bold c}_1\Vert =
\dots =\Vert{\bold b}-{\bold c}_i\Vert= R(F_{d-i}).$$ From this it is clear 
that
$$R(F_{d-i})=\text{\rm dist}({\bold o}, \text{\rm aff}F_{d-i})=\text{\rm dist}
({\bold o}, {\bold b}) \text{\rm \ with\ }{\bold b}=C_i\cap\text{\rm aff}
F_{d-i}.$$
We are left to show that ${\bold b}\in\text{\rm relint}F_{d-i}$.

\OneFigure{file=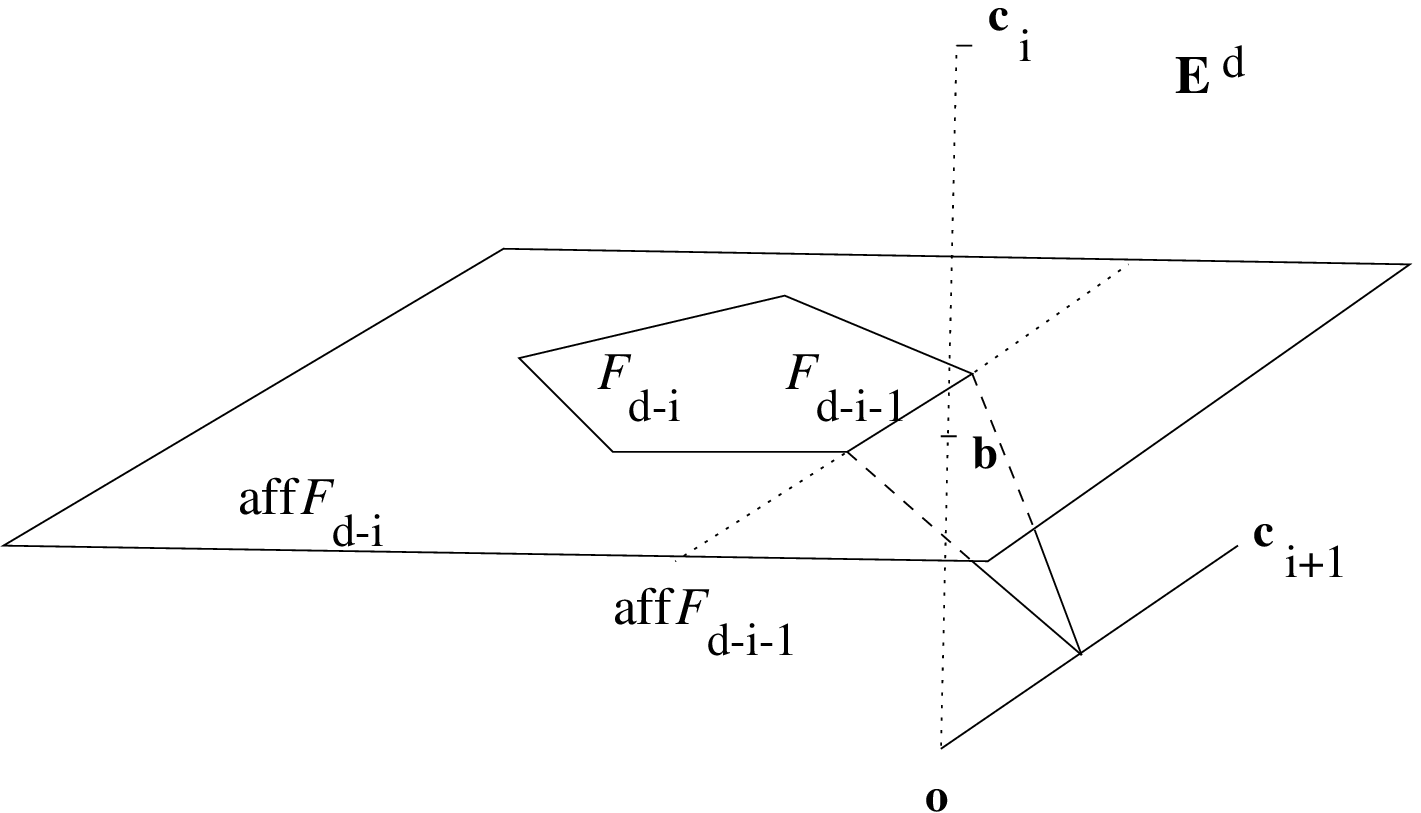}{Figure 2}

Assume that ${\bold b}\notin\text{\rm relint}F_{d-i}$ (see Fig. 2). 
Then there is a 
$(d-i-1)-$dimensional face of $F_{d-i}$ say, $F_{d-i-1}$ with the property that 
$\text{\rm aff}F_{d-i-1}$ separates the point ${\bold b}$ from
$F_{d-i}$ in $\text{\rm aff}F_{d-i}$. This means that there is a unit ball
of $\Cal P$ with center say, ${\bold c}_{i+1}$ such that the perpendicular
bisector of the line segment ${\bold o}{\bold c}_{i+1}$ intersects $\text{\rm aff}F_{d-i}$ in $\text{\rm aff}F_{d-i-1}$ and therefore it separates $F_{d-i}$
from the point ${\bold b}$ in ${\bold E^d}$.

As a result we get that
$$\Vert {\bold b}-{\bold c}_{i+1}\Vert\le \Vert {\bold b}-{\bold o}\Vert =
\Vert{\bold b}-{\bold c}_1\Vert = \dots =\Vert {\bold b}-{\bold c}_i\Vert =
R(F_{d-i}).$$ 
Now, let $H^+_k$ (resp., $H^+_0$) denote the closed halfspace of ${\bold E^d}$
that contains the point ${\bold c}_k$ (resp., ${\bold o}$) and is bounded by the hyperplane of ${\bold E^d}$ that passes through the point ${\bold b}$ and is
perpendicular to the vector ${\bold b}-{\bold c}_k$ (resp., ${\bold b}$), where
$1\le k\le i$. As ${\bold b}\in C_i$ it is easy to see that
$$H^+_0\cup H^+_1\cup \dots \cup H^+_i = {\bold E^d}.$$
From this then it is immediate that there exists a $k$ with $0\le k\le i$ such
that ${\bold c}_{i+1}\in H^+_k$. In other words, there is a point 
${\bold c}\in\{ {\bold o}, {\bold c}_1, \dots , {\bold c}_i\}$ such that
$$\angle{\bold c}{\bold b}{\bold c}_{i+1}\le \frac{\pi}{2}.$$
As $\Vert {\bold b}-{\bold c}_{i+1}\Vert\le \Vert {\bold b}-{\bold c}\Vert =
R(F_{d-i})<\sqrt{2}$ this implies in a straightforward way that
$\Vert {\bold c}-{\bold c}_{i+1}\Vert < 2$, a contradiction.
This completes the proof of Lemma 2. \qed

\heading 2. Wedges of type I, II, III, \\ truncated wedges of type I, II, and 
            \\  some of their metric properties
\endheading

Let $F_0\subset F_1\subset \dots \subset F_{d-1}$ be an arbitrary flag of the
Voronoi polytope $P$. Then let ${\bold r}_i\in F_{d-i}$ be the uniquely determined
point of the $(d-i)-$dimensional face $F_{d-i}$ of $P$ that is closest to
the center point ${\bold o}$ of $P$ that is let
$${\bold r}_i\in F_{d-i}\text{\rm \ such \ that\ }
\Vert{\bold r}_i\Vert =\text{\rm min}\{\Vert{\bold x}\Vert\ \vert\ {\bold x}\in
F_{d-i}\}, \text{\rm \ where \ }1\le i\le d.$$
\subhead Definition 1
\endsubhead
If the vectors ${\bold r}_1, \dots , {\bold r}_i$ are linearly independent
in $\bold E^d$, then
we call 
$\text{\rm conv}\{{\bold o}, {\bold r}_1, \dots , {\bold r}_i\}$ 
the {\it $i-$dimensional Rogers simplex} assigned to the subflag 
$F_{d-i}\subset \dots \subset F_{d-1}$ of the Voronoi polytope $P$, 
where $1\le i\le d$. If $\text{\rm conv}\{{\bold o}, {\bold r}_1, \dots ,
{\bold r}_d\}\subset {\bold E^d}$ is the $d-$dimensional Rogers simplex
assigned to the flag $F_0\subset\dots\subset F_{d-1}$ of $P$, then
$\text{\rm conv}\{{\bold r}_{d-i}, \dots , {\bold r}_d\}$ is called 
the {\it $i-$dimensional base} of the given d-dimensional Rogers simplex and
$\text{\rm dist}({\bold o}, \text{\rm aff}\{{\bold r}_{d-i}, \dots , {\bold r}_d\})=\text{\rm dist}({\bold o}, \text{\rm aff}F_i)=R(F_i)$ is called 
the {\it height} assigned to the $i-$dimensional base, where $1\le i\le d$.
\bigskip
\subhead Definition 2
\endsubhead
The $i-$dimensional simplex $Y=\text{\rm conv}\{{\bold o}, {\bold y}_1,
\dots , {\bold y}_i\}\subset {\bold E^d}$ with vertices ${\bold y}_0={\bold o}, {\bold y}_1, \dots , {\bold y}_i $ is called an 
{\it $i-$dimensional orthoscheme} if for each
$j, 0\le j\le i-1$ the vector ${\bold y}_j$ is orthogonal to the linear hull
$\text{\rm lin}\{ {\bold y}_k-{\bold y}_j\ \vert \ j+1\le k\le i\}$, where
$1\le i\le d$.
\bigskip
 It is shown in \cite{23} that the union of the $d-$dimensional 
Rogers simplices 
of the Voronoi polytope $P$ is the polytope $P$ itself and their interiors
are pairwise disjoint. This fact together with Corollary 1 and Lemma 2 imply 
the following metric properties of Rogers simplices in a 
straightforward way.
\proclaim {Lemma 3}
\roster
\item If $\text{\rm conv}\{{\bold o}, {\bold r}_1, \dots , {\bold r}_i\}$ is 
an {\it $i-$dimensional Rogers simplex} assigned to the subflag 
$F_{d-i}\subset \dots \subset F_{d-1}$ of the Voronoi polytope $P$, then
$\sqrt{\frac{2j}{j+1}}\le \Vert{\bold r}_j\Vert$ for all
$1\le j\le i$, where $1\le i\le d$. 
\item If $F_{d-i}\subset\dots\subset F_{d-1}$ is a subflag
of the Voronoi polytope $P$ with $R(F_{d-i})<\sqrt{2}$, then
$\text{\rm conv}\{{\bold o}, {\bold r}_1, \dots , {\bold r}_i\}$ is
an $i-$dimensional Rogers simplex which is in fact, an $i-$dimensional
orthoscheme (in short, an $i-$dimensional Rogers orthoscheme) with the property that each
${\bold r}_j\in\text{\rm relint}F_{d-j}, 1\le j\le i$ is the orthogonal projection
of $\bold o$ onto $\text{\rm aff}F_{d-j}$, where 
$1\le i\le d$.
\item If $F_{2}\subset\dots\subset F_{d-1}$ is a subflag
of the Voronoi polytope $P\subset {\bold E^d}, 3\le d$ with $R(F_{2})<\sqrt{2}$,
then the union of the $2-$dimensional bases of the $d-$dimensional Rogers
simplices that contain the orthoscheme $\text{\rm conv}\{{\bold o},
{\bold r}_1,$ $\dots , {\bold r}_{d-2}\}$ is the (uniquely determined) $2-$dimensional face $F_2$ of the Vo\-ro\-no\-i polytope $P$ that is totally orthogonal to $\text{\rm conv}\{{\bold o}, {\bold r}_1, \dots , {\bold r}_{d-2}\}$ at the point ${\bold r}_{d-2}$ and so, $\Vert {\bold r}_{d-2}\Vert=
\text{\rm dist}({\bold o}, \text{\rm aff}F_2)$ with ${\bold r}_{d-2}\in
\text{\rm relint}F_2$.
\endroster
\endproclaim
Now, we are ready for the definitions of wedges and truncated wedges. (As an illustration see Fig. 3.) Recall that for any $2-$dimensional face $F_2$ of the Voronoi polytope $P\subset{\bold E^d}, 3\le d$ we have that $\sqrt{\frac{2(d-2)}{d-1}}\le R(F_2)$.
\bigskip
\subhead Definition 3
\endsubhead
\roster
\item Let $F_2$ be a $2-$dimensional face of the Voronoi polytope $P\subset
{\bold E^d}, 3\le d$ with
$\sqrt{\frac{2(d-2)}{d-1}}\le R(F_2)<\sqrt{\frac{2(d-1)}{d}}$ and let
$\text{\rm conv}\{{\bold o}, {\bold r}_1, \dots , {\bold r}_{d-2}\}$ be any
$(d-2)-$dimensional Rogers simplex  
with ${\bold r}_{d-2}\in \text{\rm relint}F_{2}$. Then the union $W_I$ of the $d-$dimensional Rogers simplices of $P$ that contain the orthoscheme 
$\text{\rm conv}\{{\bold o}, {\bold r}_1, \dots , {\bold r}_{d-2}\}$ is 
called a {\it wedge of type I (generated by the $(d-2)-$dimensional Rogers orthoscheme $\text{\rm conv}\{{\bold o}, {\bold r}_1, \dots , {\bold r}_{d-2}\}$)}. $F_2$ is called the {\it $2-$dimensional base} of $W_I$ and $\Vert {\bold r}_{d-2}\Vert=
\text{\rm dist}({\bold o}, \text{\rm aff}F_2)$ is the {\it height} of $W_I$ assigned to the base $F_2$.
\item Let $F_2$ be a $2-$dimensional face of the Voronoi polytope $P\subset
{\bold E^d}, 3\le d$ with
$\sqrt{\frac{2(d-1)}{d}}\le R(F_2)<\sqrt{\frac{2d}{d+1}}$ and let
$\text{\rm conv}\{{\bold o}, {\bold r}_1, \dots , {\bold r}_{d-2}\}$ be any
$(d-2)-$dimensional Rogers simplex  
with ${\bold r}_{d-2}\in \text{\rm relint}F_{2}$. Then the union $W_{II}$ of the $d-$dimensional Rogers simplices of $P$ that contain the orthoscheme 
$\text{\rm conv}\{{\bold o}, {\bold r}_1, \dots , {\bold r}_{d-2}\}$ is 
called a {\it wedge of type II (generated by the $(d-2)-$di\-men\-si\-o\-nal Rogers
orthoscheme $\text{\rm conv}\{{\bold o}, {\bold r}_1, \dots , {\bold r}_{d-2}\}$)}. $F_2$ is called the {\it $2-$di\-men\-si\-o\-nal base} of $W_{II}$ and $\Vert {\bold r}_{d-2}\Vert=
\text{\rm dist}({\bold o}, \text{\rm aff}F_2)$ is the {\it height} of $W_{II}$ assigned to the base $F_2$.
\item Let $\text{\rm conv}\{{\bold o}, {\bold r}_1, \dots , {\bold r}_d\}$ be
the $d-$dimensional Rogers simplex assigned to the flag $F_0\subset F_1\dots\subset F_{d-1}$ of the Voronoi polytope $P\subset
{\bold E^d}, 3\le d$ with $\sqrt{\frac{2d}{d+1}}\le R(F_2)$. Then
$W_{III}=\text{\rm conv}\{{\bold o}, {\bold r}_1, \dots , {\bold r}_d\}$ is called a {\it wedge of type III}. 

\endroster
\bigskip
It useful to recall that for any vertex $F_0$ of the Voronoi polytope $P
\subset {\bold E^d}$ we have that $\sqrt{\frac{2d}{d+1}}\le R(F_0)$.

\subhead Definition 4
\endsubhead
Let $\overline{B}=\{ {\bold x}\in {\bold E^d}\vert \ \text{dist}({\bold o}, {\bold x})=\Vert{\bold x}\Vert\le \sqrt{\frac{2d}{d+1}}\}$.
\roster
\item If $W_I$ is a wedge of type I with the $2-$dimensional base $F_2$ which is generated by the $(d-2)-$dimensional Rogers 
orthoscheme $\text{\rm conv}\{{\bold o}, {\bold r}_1, \dots , 
{\bold r}_{d-2}\}$ of the Voronoi polytope $P\subset {\bold E^d}, 3\le d$, then 
$$\overline{W}_I=\text{\rm conv}\left((\overline{B}\cap F_2)\cup\{{\bold o}={\bold r}_0, \dots , {\bold r}_{d-3}\}\right)$$ 
is called the {\it truncated wedge of type I} with the
{\it $2$-dimensional base} $\overline{B}\cap F_2$ generated by the $(d-2)-$dimensional Rogers orthoscheme 
$$\text{\rm conv}\{{\bold o}, {\bold r}_1, \dots , {\bold r}_{d-2}\}.$$

\item If $W_{II}$ is a wedge of type II with the $2-$dimensional base $F_2$ which is generated by the $(d-2)-$dimensional Rogers 
orthoscheme $\text{\rm conv}\{{\bold o}, {\bold r}_1, \dots , 
{\bold r}_{d-2}\}$ of the Voronoi polytope $P\subset {\bold E^d}, 3\le d$, 
then $$\overline{W}_{II}=\text{\rm conv}\left((\overline{B}\cap F_2)\cup\{{\bold o}={\bold r}_0, \dots , {\bold r}_{d-3}\}\right) $$ 
is called the {\it truncated wedge of type II}
with the {\it $2$-dimensional base} $\overline{B}\cap F_2$ generated by the $(d-2)-$dimensional Rogers orthoscheme 
$$\text{\rm conv}\{{\bold o}, {\bold r}_1, \dots , {\bold r}_{d-2}\}.$$
\endroster
\bigskip

As the following claim can be proved by Lemma 3 in a straightforward way, we omit its simple proof (see also Fig. 3).
\proclaim {Sublemma 2}
\roster
\item Let $W_I$ (resp., $W_{II}$) denote the wedge of type I (resp., of type II)
with the $2-$dimensional base $F_2$ which is generated by the $(d-2)-$dimensional Rogers 
orthoscheme $\text{\rm conv}\{{\bold o}, {\bold r}_1, \dots , 
{\bold r}_{d-2}\}$ of the Voronoi polytope $P\subset {\bold E^d}, 3\le d$.
If the points ${\bold x}, {\bold y}\in \text{\rm aff}F_2$ are choosen so that
the triangle $\triangle{\bold r}_{d-2}{\bold x}{\bold y}$ has a right angle
at the vertex ${\bold x}$, then $\text{\rm conv}\{{\bold o}, {\bold r}_1,
\dots ,{\bold r}_{d-2}, {\bold x}, {\bold y}\}$ is a $d-$di\-men\-si\-o\-nal
orthoscheme. Moreover, if ${\bold z}\in \text{\rm aff}F_2$ is an arbitrary point,
then $\text{\rm conv}\{{\bold o}={\bold r}_0, \dots , {\bold r}_{d-3}, {\bold z}\}$
is a $(d-2)-$dimensional orthoscheme.
\item Let $W_I$ denote the wedge of type I
with the $2-$dimensional base $F_2$ which is generated by the $(d-2)-$dimensional Rogers orthoscheme $\text{\rm conv}\{{\bold o}={\bold r}_0, {\bold r}_1, \dots ,$ 
${\bold r}_{d-2}\}$ of the Voronoi polytope $P\subset {\bold E^d}, 3\le d$.
Let $Q_2\subset{\text{\rm aff}}F_2$ and $Q^*_2\subset{\text{\rm aff}}F_2$
be compact convex sets with $\text{\rm relint}Q_2\cap\text{\rm relint}Q^*_2=\emptyset$. If
$K_2=Q_2$ (resp., $K^*_2=Q^*_2$) and $K_j=\text{\rm conv}(K_{j-1}\cup\{{\bold r}_{d-j}\})$
(resp., $K^*_j=\text{\rm conv}(K^*_{j-1}\cup\{{\bold r}_{d-j}\})$) for $j=3, \dots , d$, then
$K_d=\text{\rm conv}(Q_2\cup\{{\bold o}={\bold r}_0, \dots , {\bold r}_{d-3}\})$
(resp., $K^*_d=\text{\rm conv}(Q^*_2\cup\{{\bold o}={\bold r}_0, \dots , {\bold r}_{d-3}\})$)
moreover, $\text{\rm relint}K_d\cap\text{\rm relint}K^*_d=\emptyset$. A similar statement holds
for $W_{II}$.
 
\item Let $W_I$ (resp., $\overline{W}_I$) denote the wedge of type I (resp.,
truncated wedge of type I) with the $2-$dimensional base $F_2$ (resp., $\overline{B}\cap F_2$) 
which is generated by the $(d-2)-$dimensional Rogers orthoscheme $\text{\rm conv}\{{\bold o}={\bold r}_0, {\bold r}_1, \dots , {\bold r}_{d-2}\}$ of the Voronoi polytope $P\subset {\bold E^d}, 3\le d$. If $K_2=F_2$ 
(resp., $K_2=\overline{B}\cap F_2$) and $K_j=\text{\rm
conv} (K_{j-1}\cup \{{\bold r}_{d-j}\})$ for $j=3, \dots , d$, then
$K_d=W_I$ (resp., $K_d=\overline{W}_I$). Similar statements hold for $W_{II}$ and $\overline{W}_{II}$. 

\endroster
\endproclaim

\OneFigure{file=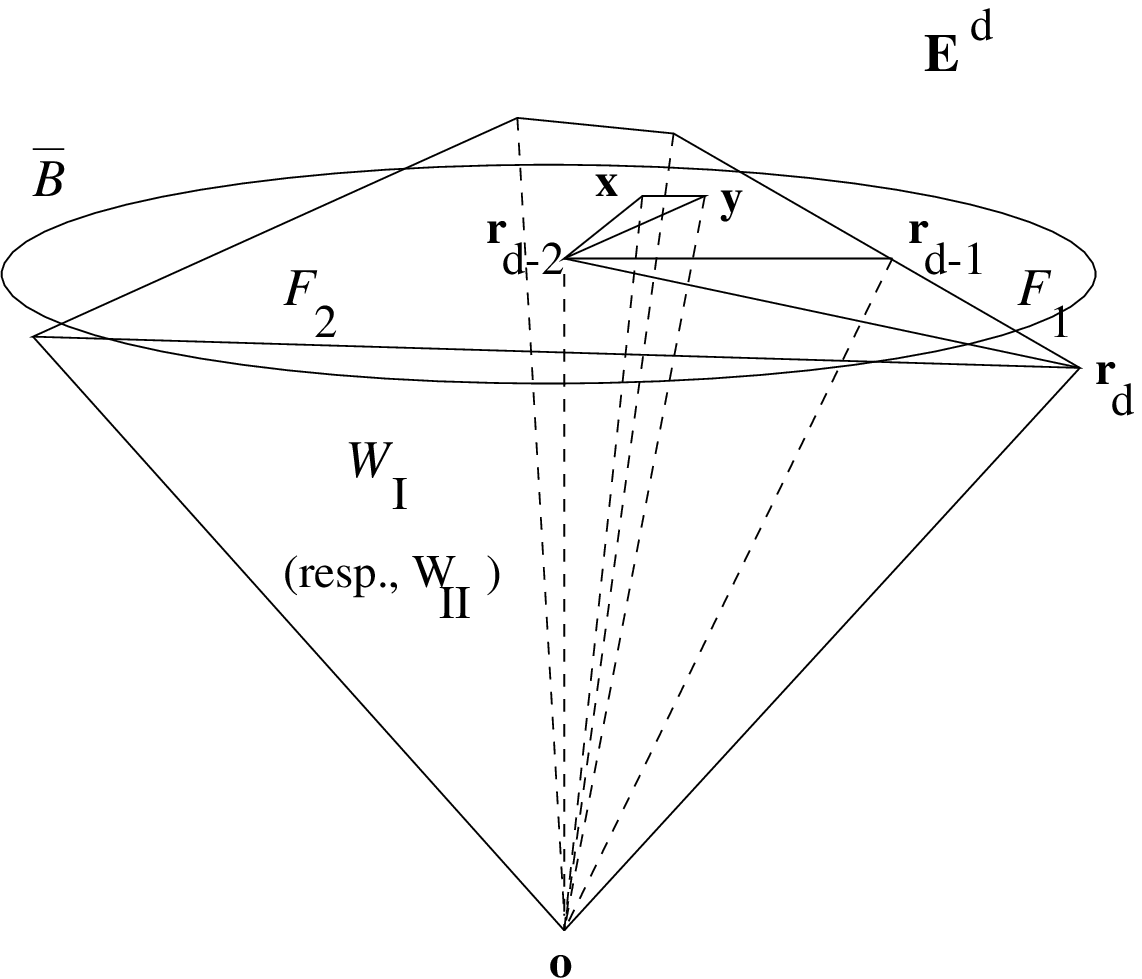}{Figure 3}

The core part of this section is Corollary 2 that follows from Lemma 4
in a trivial way.
\proclaim {Lemma 4}
Let $F_2$ be an arbitrary $2-$dimensional face of the Voronoi polytope $P\subset
{\bold E^d}$ of dimension $d\ge 8$. Then the number of sides $F_1$ of the face $F_2$ with $R(F_1)<\sqrt{\frac{2d}{d+1}}$ is at most $4$.
\endproclaim
\subhead Proof
\endsubhead
Assume that there are $5$ sides (i.e. boundary line segments) say, $E_1, E_2$, $E_3$, $E_4$ and $E_5$ of the face $F_2$ of the Voronoi polytope $P\subset {\bold E^d}, d\ge 8$ with the property that 
$$R(E_m)=\text{\rm dist}({\bold o}, \text{\rm aff}E_m) <\sqrt{\frac{2d}{d+1}}, \text{\rm \ for\  all \ }
1\le m\le 5.$$

First, (as we have seen above) we can pick $(d-2)$ unit balls of $\Cal P$ with centers say, ${\bold c}_1,
\dots , {\bold c}_{d-2}$ such that each of their Voronoi cells contains 
$F_2$ (and so $\text{\rm aff}\{{\bold o}, {\bold c}_1, \dots , {\bold c}_{d-2}\}=\text{\rm lin}
\{{\bold c}_1, \dots , {\bold c}_{d-2}\}$ is totally orthogonal to $\text{\rm aff}F_2$) moreover,
the point ${\bold b}=\text{\rm lin}\{{\bold c}_1, \dots , {\bold c}_{d-2}\}\cap \text{\rm aff}F_2$
belongs to the $(d-2)-$dimensional simplex $C_{d-2}=\text{\rm conv}\{{\bold o}, {\bold c}_1, \dots , 
{\bold c}_{d-2}\}$. (Also, notice that ${\bold b}$ is the center of the $(d-3)-$dimensional sphere
that passes through the vertices of $C_{d-2}$.) The existence of the sides $E_m, 1\le m\le 5$ implies that $R(F_2)<\sqrt{\frac{2d}{d+1}}<\sqrt{2}$ and so, via Lemma 2 we get that ${\bold b}\in \text{\rm relint}F_2$.

Second, for each side $E_m, 1\le m\le 5$ there exists a unit ball of $\Cal P$ with center say, ${\bold c}_{f(m)}$ such that its Voronoi 
cell intersects $F_2$ in $E_m$. Let the hyperplane of ${\bold E^d}$
spanned by the $(d-2)$-dimensional simplex $C_{d-2}$ and the center point ${\bold c}_{f(m)}$ be denoted
by $H_{f(m)}$. As the hyperplanes $H_{f(m)}, 1\le m\le 5$ all contain  the $(d-2)-$dimensional linear subspace $\text{\rm lin}\{{\bold c}_1, \dots , {\bold c}_{d-2}\}$ there must be $i$ and $j$ such that
the angle between $H_i$ and $H_j$ is $\frac {2\pi}{5}$ or less. We complete the proof by showing that this forces $\text{\rm dist}({\bold c}_i, {\bold c}_j)<2$, an impossibility.

Let ${\bold a}_i$ (resp., ${\bold a}_j$) be the center of the $(d-2)-$dimensional sphere that passes through the points ${\bold o}, 
{\bold c}_1, \dots , {\bold c}_{d-2}, {\bold c}_i$ (resp., ${\bold o}, 
{\bold c}_1, \dots , {\bold c}_{d-2}, {\bold c}_j$) (see Fig. 4). By assumption,
$\Vert {\bold a}_i\Vert=\text{\rm dist}({\bold o}, {\bold a}_i)=R(E_{f^{-1}(i)})<\sqrt{\frac{2d}{d+1}}$ (resp., $\Vert {\bold a}_j\Vert=\text{\rm dist}({\bold o}, {\bold a}_j)=R(E_{f^{-1}(j)})<\sqrt{\frac{2d}{d+1}}$) and so Lemma 2 implies that
${\bold a}_i\in\text{\rm relint}E_{f^{-1}(i)}$ (resp., ${\bold a}_j\in\text{\rm relint}E_{f^{-1}(j)}$). Now, let ${\bold a}_i^*$ (resp., ${\bold a}_j^*$) be
the uniquely determined point on the halfline that emanates from ${\bold b}$
and passes through ${\bold a}_i$ (resp., ${\bold a}_j$) with the property that
$\Vert{\bold a}_i^*\Vert =\sqrt{\frac{2d}{d+1}}$ (resp., $\Vert{\bold a}_j^*\Vert =\sqrt{\frac{2d}{d+1}}$). Moreover, if $S_i^{d-2}=
\{{\bold x}\in
H_i\ \vert \ \Vert {\bold x}-{\bold a}_i^*\Vert=\sqrt{\frac{2d}{d+1}}
\}$ (resp., $S_j^{d-2}=\{{\bold x}\in
H_j\ \vert \ \Vert {\bold x}-{\bold a}_j^*\Vert=\sqrt{\frac{2d}{d+1}}
\}$), then let ${\bold c}_i^*$ (resp., ${\bold c}_j^*$) be the intersection
of $S_i^{d-2}$ (resp., $S_j^{d-2}$) with the halfline emanating from ${\bold b}$
and passing through the point ${\bold c}_i$ (resp., ${\bold c}_j$). Finally,
let ${\bold b}_i$ (resp., ${\bold b}_j$) be the orthogonal projection of the
point ${\bold c}_i^*$ (resp., ${\bold c}_j^*$) onto the plane $\text{\rm aff}F_2$. As $\text{\rm lin}\{{\bold c}_1, \dots ,{\bold c}_{d-2}\}$ is orthogonal to $\text{\rm aff}F_2$ and  $\text{\rm lin}\{{\bold c}_1, \dots ,{\bold c}_{d-2}\}\subset H_i$ (resp., $\text{\rm lin}\{{\bold c}_1, \dots ,{\bold c}_{d-2}\}\subset H_j$) the points ${\bold b}, {\bold a}_i, 
{\bold a}_i^*, {\bold b}_i$ (resp., ${\bold b}, {\bold a}_j, 
{\bold a}_j^*, {\bold b}_j$) lie on a line of $\text{\rm aff}F_2$.

\OneFigure{file=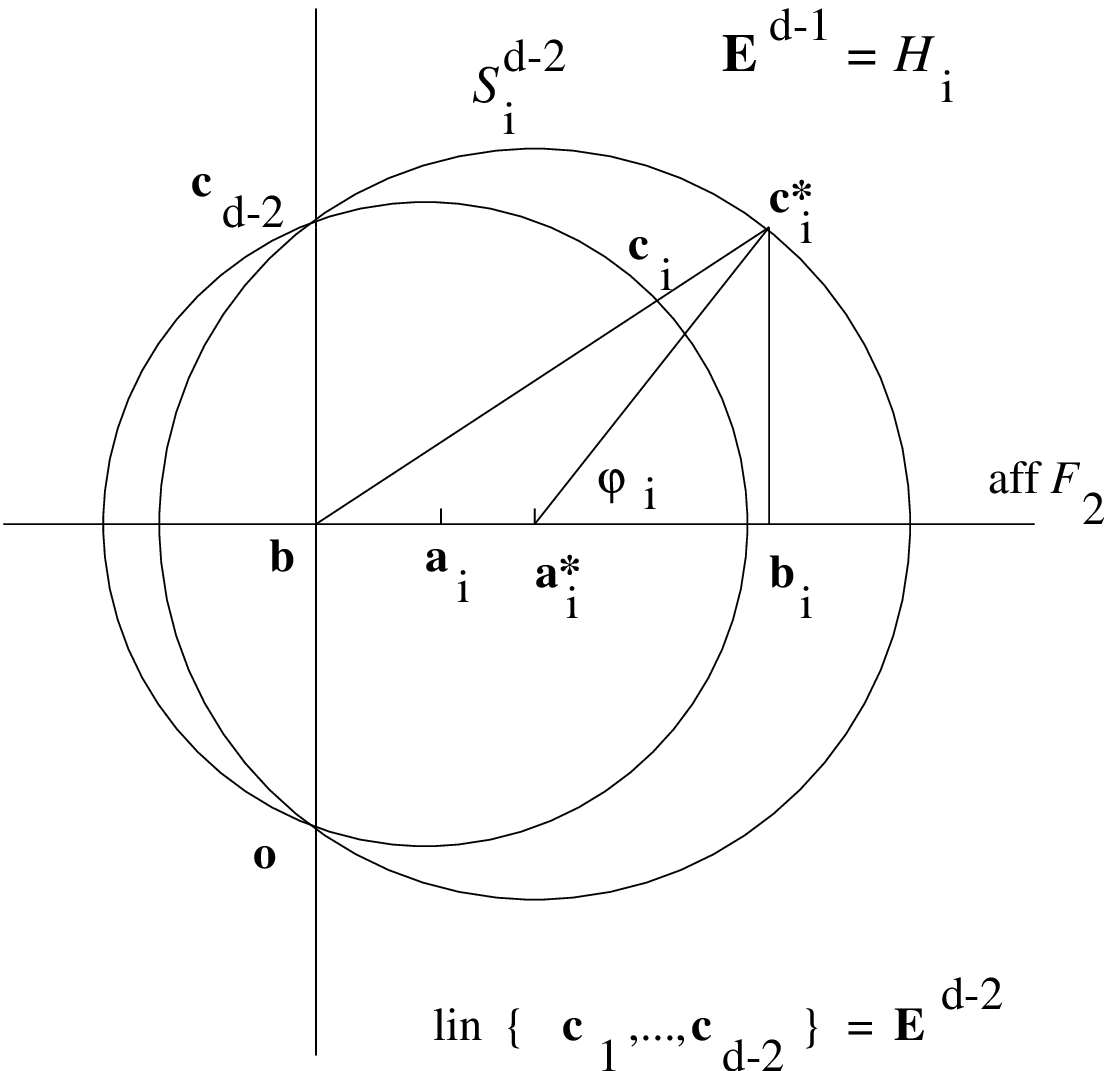}{Figure 4}

\proclaim {Sublemma 3}
$\text{\rm dist}({\bold c}_i, {\bold c}_j)< \text{\rm dist}({\bold c}_i^*, {\bold c}_j^*).$
\endproclaim

\subhead Proof
\endsubhead
Let $s_i=\Vert{\bold c}_i-{\bold b}\Vert , s_i^*=\Vert{\bold c}_i^*-{\bold b}
\Vert$ (resp.,  $s_j=\Vert{\bold c}_j-{\bold b}\Vert , s_j^*=\Vert{\bold c}_j^*-{\bold b}\Vert$). If $\psi=\angle{\bold c}_i{\bold b}{\bold c}_j=
\angle{\bold c}_i^*{\bold b}{\bold c}_j^*$, then
$$\Vert{\bold c}_i-{\bold c}_j\Vert^2=s_i^2+s_j^2-2s_is_j\cos \psi
\text {\rm \ and \ }$$
$$\Vert{\bold c}_i^*-{\bold c}_j^*\Vert^2=(s_i^*)^2+(s_j^*)^2-2(s_i^*)(s_j^*)\cos \psi .$$
As $s_i < s_i^*$ and $s_j < s_j^*$ it is sufficient to show that
$$\frac{\partial{ }}{\partial{s_i}}\Vert{\bold c}_i-{\bold c}_j\Vert^2=
2(s_i-s_j\cos \psi)>0 \text{\rm \ and \ }$$
$$\frac{\partial{ }}{\partial{s_j}}\Vert{\bold c}_i-{\bold c}_j\Vert^2=
2(s_j-s_i\cos \psi)>0 .$$
By symmetry it is sufficient to show that $s_i-s_j\cos \psi>0$. Assume that
$s_i\le s_j\cos \psi$. Then on the one hand, $\frac{\pi}{2}\le \angle {\bold b}{\bold c}_i{\bold c}_j$ and so $4< s_i^2+\Vert {\bold c}_i-{\bold c}_j\Vert ^2\le s_j^2$. On the other hand, $s_j\le \Vert {\bold a}_j-{\bold b}\Vert +
\Vert {\bold c}_j - {\bold a}_j\Vert =\sqrt{\Vert {\bold a}_j\Vert^2 - \Vert
{\bold b}\Vert ^2} +\Vert {\bold a}_j\Vert$. By assumption and Corollary 1 we
know that $\Vert {\bold a}_j\Vert <\sqrt {\frac{2d}{d+1}}$ and $\sqrt{
\frac{2(d-2)}{d-1}}\le \Vert {\bold b}\Vert$ and so
$$s_j < \sqrt{\frac{2d}{d+1}-\frac{2(d-2)}{d-1} }+\sqrt{\frac{2d}{d+1}}\le
2 \text{\rm \ for \ all\ }d\ge 3,$$
a contradiction. This completes the proof of Sublemma 3. $\qed$

\bigskip

Thus, it will suffice to prove that $\text{\rm dist}({\bold c}_i^*, {\bold c}_j^*)\le 2$.

Let $\varphi _i=\pi -\angle {\bold b}{\bold a}_i^*{\bold c}_i^*$ (resp., $\varphi _j=\pi -\angle {\bold b}{\bold a}_j^*{\bold c}_j^*$). 
Sublemma 4, which we prove at the end of this section, shows that
$0\le \varphi _i\le \frac{\pi}{2}$ (resp., $0\le \varphi _j\le \frac{\pi}{2}$) as indicated in
Fig. 4.

If $\phi =\angle {\bold a}_i{\bold b}{\bold a}_j$ and $l_i=\Vert {\bold b}_i-{\bold b}\Vert$ (resp., $l_j=\Vert {\bold b}_j-{\bold b}\Vert$), then we get that
$$\Vert {\bold c}_i^*-{\bold c}_j^*\Vert ^2\le \Vert {\bold b}_i-{\bold b}_j\Vert ^2+ 
(\Vert {\bold c}_i^*-{\bold b}_i\Vert + \Vert {\bold c}_j^*-{\bold b}_j\Vert)^2=$$
$$(l_i^2+l_j^2-2l_il_j\cos \phi )+\frac{2d}{d+1}(\sin \varphi_i +\sin \varphi_j)^2.$$
By assumption $0\le \phi \le \frac{2\pi}{5}$ and so
$$\Vert {\bold c}_i^*-{\bold c}_j^*\Vert ^2\le (l_i^2+l_j^2-2l_il_j\cos \frac{2\pi}{5} ) +\frac{2d}{d+1}(\sin \varphi_i +\sin \varphi_j)^2.$$

Substituting  $l_i=\sqrt{\frac{2d}{d+1}-\Vert {\bold b}\Vert^2 }+\sqrt{\frac{2d}{d+1}}\cos \varphi_i$
(resp., $l_j=\sqrt{\frac{2d}{d+1}-\Vert {\bold b}\Vert^2 }+\sqrt{\frac{2d}{d+1}}\cos \varphi_j$)
yields
$$\Vert {\bold c}_i^*-{\bold c}_j^*\Vert ^2\le
\left(2-\cos\frac{2\pi}{5}\right)\frac{4d}{d+1} -2\left(1-\cos\frac{2\pi}{5}\right)\Vert {\bold b}\Vert^2+$$
$$\frac{4d}{d+1}\sin\varphi_i\sin\varphi_j-\left(\frac{4d}{d+1}\cos\frac{2\pi}{5}\right)\cos\varphi_i
\cos\varphi_j+$$
$$2\left(1-\cos\frac{2\pi}{5}\right)\sqrt{\frac{2d}{d+1}}
\left(\cos\varphi_i+\cos\varphi_j\right)\sqrt{\frac{2d}{d+1}-\Vert{\bold b}\Vert^2 }.$$

Corollary 1 implies that $\sqrt{\frac{2(d-2)}{d-1}}\le\Vert{\bold b}\Vert$ and so
$$\Vert {\bold c}_i^*-{\bold c}_j^*\Vert ^2\le
\left(2-\cos\frac{2\pi}{5}\right)\frac{4d}{d+1}- \left(1-\cos\frac{2\pi}{5}\right)\frac{4(d-2)}{d-1}+$$
$$\frac{4d}{d+1}\sin\varphi_i\sin\varphi_j-\left(\frac{4d}{d+1}\cos\frac{2\pi}{5}\right)\cos\varphi_i
\cos\varphi_j+$$
$$2\left(1-\cos\frac{2\pi}{5}\right)\sqrt{\frac{2d}{d+1}}
\left(\cos\varphi_i+\cos\varphi_j\right)\sqrt{\frac{2d}{d+1}-\frac{2(d-2)}{d-1} }=G(\varphi_i,\varphi_j).$$

Then straightforward computation yields 
$$\frac{\partial { }}{\partial{\varphi_i}}G(\varphi_i, \varphi_j)=
\frac{4d}{d+1}\sin\varphi_j\cos\varphi_i +$$ $$\sin\varphi_i\left(\frac{4d}{d+1}\cos\frac{2\pi}{5}\cos\varphi_j
-4\left(1-\cos\frac{2\pi}{5}\right)\sqrt{\frac{2d}{(d+1)^2(d-1)} }\right)\text{\rm and}$$

$$\frac{\partial { }}{\partial{\varphi_j}}G(\varphi_i, \varphi_j)=
\frac{4d}{d+1}\sin\varphi_i\cos\varphi_j +$$ $$\sin\varphi_j\left(\frac{4d}{d+1}\cos\frac{2\pi}{5}\cos\varphi_i
-4\left(1-\cos\frac{2\pi}{5}\right)\sqrt{\frac{2d}{(d+1)^2(d-1)} }\right).$$

Now, we complete the proof of Lemma 4 by Sublemma 4, which we are left to prove at the end of this
section. The details are as follows.

According to Sublemma 4 the angle $\varphi_i$ (resp., $\varphi_j$) is maximized when $\Vert{\bold b}\Vert=
\sqrt{\frac{2(d-2)}{d-1}}$ (i.e. when the norm of ${\bold b}$ is minimized i.e. when the $(d-2)-$dimensional simplex $C_{d-2}=\text{\rm conv}\{{\bold o}, {\bold c}_1, \dots , {\bold c}_{d-2}\}$
is a regular simplex of edge length $2$) and the unit ball
centered at ${\bold c}_i^*$ (resp., ${\bold c}_j^*$) is tangent to exactly $(d-2)$ unit balls
centered at $(d-2)$ points out of ${\bold o}, {\bold c}_1, \dots , {\bold c}_{d-2}$. In other
words, as it is claimed below in Sublemma 4 
$$\frac{\sqrt{2}}{3}\frac{2d-1}{\sqrt{d(d-1)}}\le \cos\varphi_i \text{\rm \ and \ }
\frac{\sqrt{2}}{3}\frac{2d-1}{\sqrt{d(d-1)}}\le \cos\varphi_j \text{\rm for\ all\ }d\ge 2.$$
From this it follows in a straightforward way that
$$\frac{\partial { }}{\partial{\varphi_i}}G(\varphi_i, \varphi_j)\ge 0 \text{\rm \ and \ }
\frac{\partial { }}{\partial{\varphi_j}}G(\varphi_i, \varphi_j)\ge 0 \text {\rm \ for\ all\ }
d\ge 4.$$
Therefore $G(\varphi_i, \varphi_j)$ is maximized when $\varphi_i$ and $\varphi_j$ are as
large as possible that is using Sublemma 4 again when $\cos\varphi_i =\cos\varphi_j =
\frac{\sqrt{2}}{3}\frac{2d-1}{\sqrt{d(d-1)}}$. This yields
$$\Vert {\bold c}_i^*-{\bold c}_j^*\Vert ^2\le$$
$$\frac{\left(40-32\cos\frac{2\pi}{5}\right)d^2+\left(56-64\cos\frac{2\pi}{5}\right)d+
\left(16-32\cos\frac{2\pi}{5}\right) } {9\left(d^2-1\right)}\le 4 \text{\rm \ for\ all  \ }
d\ge 8$$
as desired.  

We are left to prove the following statement.
\proclaim {Sublemma 4}
 $$\frac{\sqrt{2}}{3}\frac{2d-1}{\sqrt{d(d-1)}}\le \cos\varphi_i \text{\rm \ and \ }
\frac{\sqrt{2}}{3}\frac{2d-1}{\sqrt{d(d-1)}}\le \cos\varphi_j \text{\rm for\ all\ }d\ge 4.$$
\endproclaim
\subhead Proof
\endsubhead   The extremal problem whose solution leads to the stated inequalities is the following
one.
Let the $(d-1)-$dimensional simplex $X=\text{\rm conv}\{{\bold x}_0, {\bold x}_1, \dots , {\bold x}_{d-2},$
${\bold x}_{d-1}\}$ of all edge lengths being at least $2$ be inscribed in the $(d-2)-$dimensional
sphere $S^{d-2}=\{{\bold y}\in {\bold E^{d-1}}\ \vert\ \Vert {\bold y} - {\bold a}\Vert =
\sqrt{\frac{2d}{d+1}}\}$ centered at the point ${\bold a}$ in ${\bold E^{d-1}}, d\ge 4$. Assume that
the center ${\bold b}$ of the $(d-3)-$dimensional sphere that passes through the vertices 
${\bold x}_0, {\bold x}_1, \dots , {\bold x}_{d-2}$ of $X$ 
belongs to $\text{\rm conv}\{{\bold x}_0, {\bold x}_1, \dots , {\bold x}_{d-2}\}$. Now, if $\varphi =
\pi - \angle {\bold b}{\bold a}{\bold x}_{d-1}$, then find the maximum of $\varphi$ for the simplices
$X$. (Notice that the given metric conditions force the points ${\bold a}$ and ${\bold x}_{d-1}$ lie
in the same open halfspace bounded by the hyperplane $\text{\rm aff}\{{\bold x}_0, {\bold x}_1,
\dots , {\bold x}_{d-2}\}$ in ${\bold E^{d-1}}$.)

For the sake of simplicity let $X$ itself denote an extremal simplex i.e. a $(d-1)-$dimensional
simplex of the above type for which $\varphi$ is maximal. Now, it is easy to show -  by contradiction
and by moving the vertex ${\bold x}_{d-1}$ in the proper direction on $S^{d-2}$ - that there exists a $(d-3)-$dimensional
face of the $(d-2)-$dimensional simplex $\text{\rm conv}\{{\bold x}_0, {\bold x}_1, \dots , 
{\bold x}_{d-2}\}$ say, $\text{\rm conv}\{{\bold x}_1, \dots , {\bold x}_{d-2}\}$ such that all of
its vertices lie at distance $2$ from the point ${\bold x}_{d-1}$ i.e. 
$$\Vert {\bold x}_{1} - {\bold x}_{d-1}\Vert = \Vert {\bold x}_{2} - {\bold x}_{d-1}\Vert = \dots =
\Vert {\bold x}_{d-2} - {\bold x}_{d-1}\Vert =2.$$
But, then again it is easy to prove - by contradiction and by moving the vertex ${\bold x}_0$
in the proper direction on $S^{d-2}$ - that the vertices
of the $(d-3)-$dimensional simplex $\text{\rm conv}\{{\bold x}_1, \dots , {\bold x}_{d-2}\}$ must lie
at distance $2$ from the vertex ${\bold x}_0$ as well i.e.
$$\Vert{\bold x}_1-{\bold x}_0\Vert =\Vert {\bold x}_2-{\bold x}_0\Vert =\dots =\Vert {\bold x}_{d-2}-{\bold x}_0\Vert=2.$$ 

Now, let ${\bold c}$ be the center of the $(d-4)-$dimensional sphere $S^{d-4}$ that passes through the vertices of the $(d-3)-$dimensional simplex $\text{\rm conv}\{{\bold x}_1, \dots, {\bold x}_{d-2}\}$
in $\text{\rm aff}\{{\bold x}_1, \dots, {\bold x}_{d-2}\}$ (see Fig. 5).

\OneFigure{file=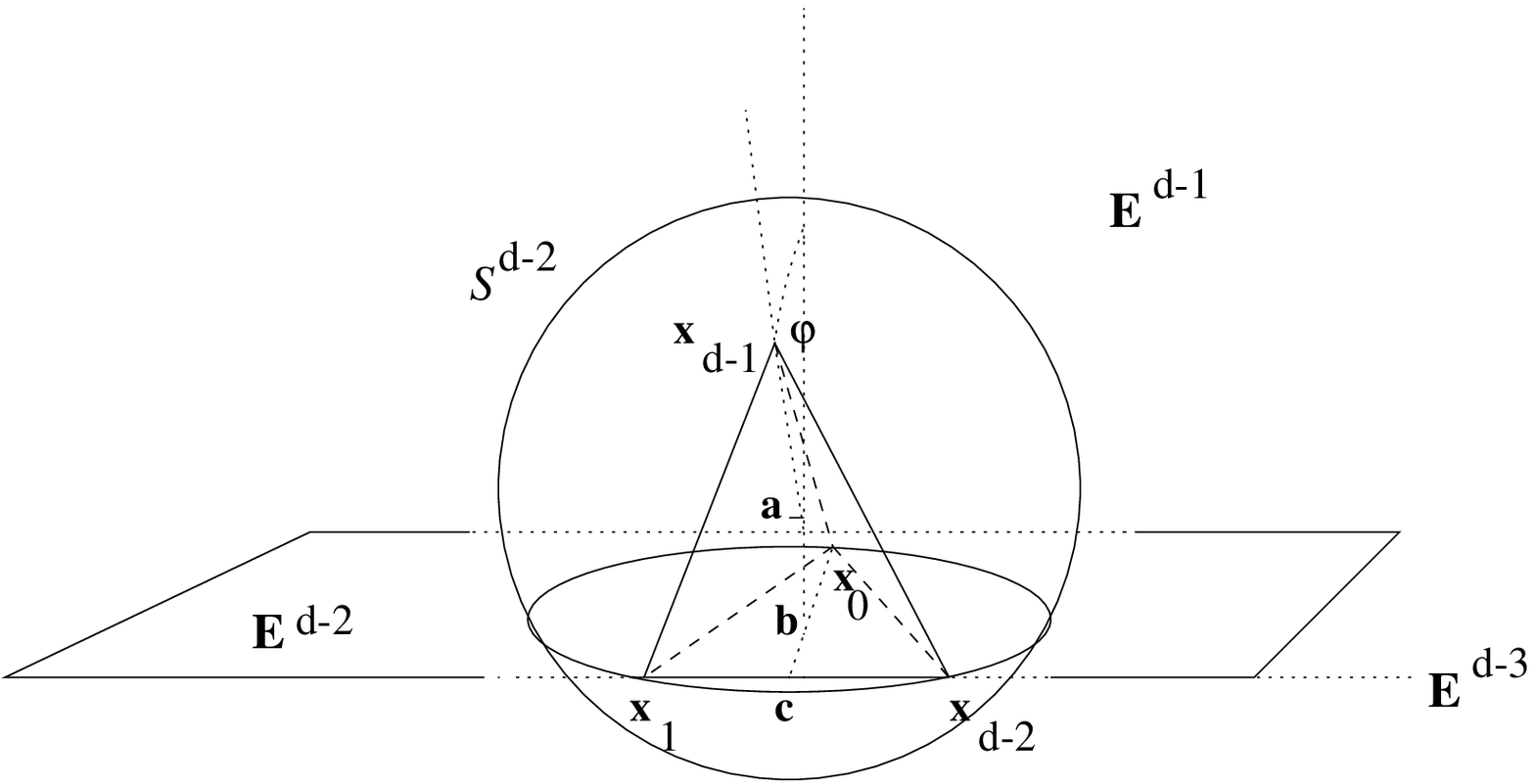}{Figure 5}

As the vertices ${\bold x}_1, \dots ,
{\bold x}_{d-2}$ lie at distance $2$ from the vertices ${\bold x}_0$ and 
${\bold x}_{d-1}$ it is clear that $\text{\rm aff}\{{\bold x}_0, {\bold x}_{d-1}, {\bold c}\}$ is
orthogonal to $\text{\rm aff}\{{\bold x}_1, \dots , {\bold x}_{d-2}\}$ (and of course, ${\bold c}=
\text{\rm aff}\{{\bold x}_0, {\bold x}_{d-1}, {\bold c}\}\cap \text{\rm aff}\{{\bold x}_1, \dots , 
{\bold x}_{d-2}\}$). Moreover, as ${\bold b}\in\text{\rm conv}\{{\bold x}_0, {\bold x}_1, \dots , {\bold x}_{d-2}\}$ and ${\bold x}_0$ (resp., ${\bold b}$) lies equidistant from the vertices ${\bold x}_1,
\dots , {\bold x}_{d-2}$ therefore ${\bold b}$ belongs to the line segment ${\bold x}_0{\bold c}$.
Finally, as ${\bold a}-{\bold b}$ is orthogonal to $\text{\rm aff}\{{\bold x}_0, {\bold x}_1,
\dots , {\bold x}_{d-2}\}$ we get that ${\bold a}\in \text{\rm aff}\{{\bold x}_0, {\bold x}_{d-1},
{\bold c}\}$ and so $\angle{\bold b}{\bold a}{\bold x}_{d-1}$ and $\varphi $ are both represented in
the plane $\text{\rm aff} \{{\bold x}_0, {\bold x}_{d-1}, {\bold c}\}$. Now, let ${\bold u}, {\bold v}$ be any pair of diametrically opposite points of $S^{d-4}$. Then by construction $\text{\rm aff}\{
{\bold u}, {\bold v}\}$ is orthogonal to $\text{\rm aff} \{{\bold x}_0, {\bold x}_{d-1}, {\bold c}\}$
and
$$\Vert {\bold u}-{\bold x}_0\Vert=\Vert {\bold u}-{\bold x}_{d-1}\Vert=
\Vert {\bold v}-{\bold x}_0\Vert=\Vert {\bold v}-{\bold x}_{d-1}\Vert=2,$$
$$\Vert {\bold u}-{\bold a}\Vert = \Vert {\bold v}-{\bold a}\Vert =
 \Vert {\bold x}_0 - {\bold a}\Vert =\Vert {\bold x}_{d-1}-{\bold a}\Vert=\sqrt{\frac{2d}{d+1}},$$
$$\Vert {\bold u}-{\bold b}\Vert = \Vert {\bold v}- {\bold b}\Vert = \Vert {\bold x}_0 - {\bold b}\Vert ,$$
moreover, Corollary 1 applied to the $(d-3)-$dimensional simplex $\text{\rm conv}\{{\bold x}_1,
\dots ,$ ${\bold x}_{d-2}\}$ yields
$$\sqrt{ \frac{2(d-3)}{d-2}}\le \Vert {\bold u}-{\bold c}\Vert=\Vert {\bold v}-{\bold c}\Vert
<\sqrt{\frac{2d}{d+1}}.$$

From now on, we work in the $3-$dimensional Euclidean space $\text{\rm aff}\{{\bold a}, {\bold b},
{\bold c}, {\bold x}_0, {\bold x}_{d-1},$ ${\bold u}, {\bold v}\}$ (see Fig. 6). Let 
$$x= \Vert {\bold u}-{\bold c}\Vert=\Vert {\bold v}-{\bold c}\Vert,$$
$$l= \Vert {\bold u}-{\bold a}\Vert = \Vert {\bold v}-{\bold a}\Vert =
 \Vert {\bold x}_0 - {\bold a}\Vert =\Vert {\bold x}_{d-1}-{\bold a}\Vert \text{\rm \ and \ }$$
$$\rho =\angle {\bold b}{\bold a}{\bold c}, \tau =\angle{\bold c}{\bold a}{\bold x}_{d-1}.$$

\centerline { }

\OneFigure{file=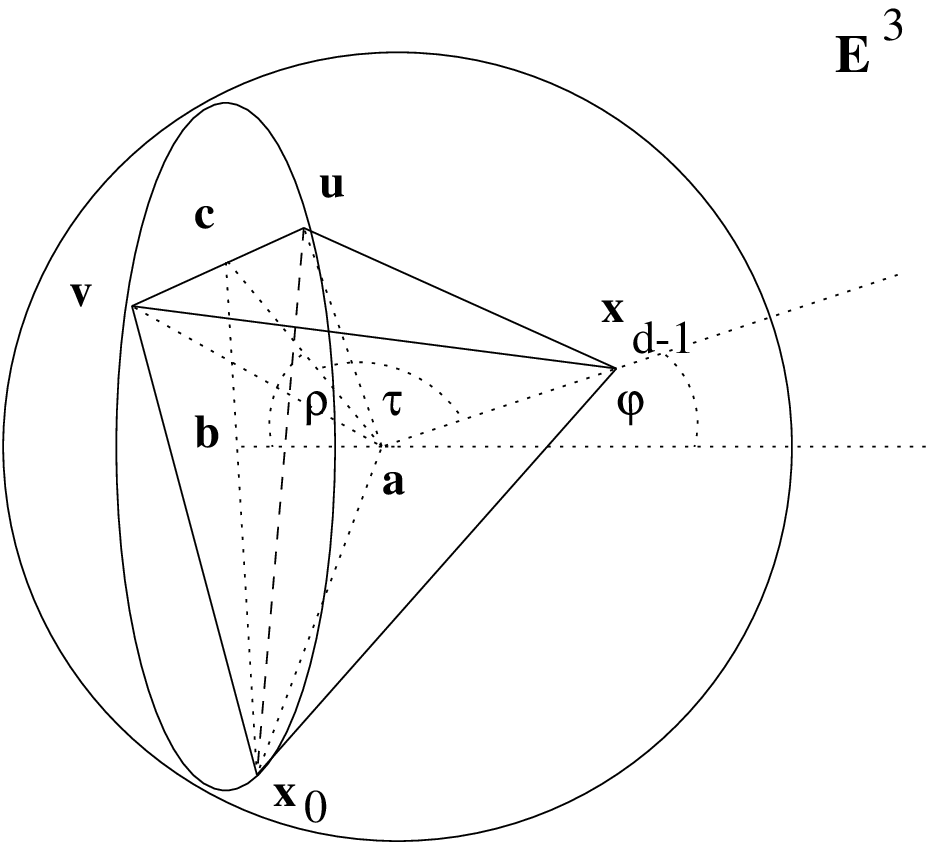}{Figure 6}

An easy elementary geometry yields
$$\cos \rho =\sqrt{\frac{l^2(4-x^2)-4}{(4-x^2)(l^2-x^2)}} \text{\rm \ and \ }
\cos \tau =\frac{l^2-2}{l\sqrt{l^2-x^2}}.$$
Therefore differentiating $\rho$ and $\tau$ as functions of $x$ we get that
$$\deriv{\rho} =\frac{-x\left[l^2x^4-8\left(l^2-1\right)x^2+12l^2-16\right] }
                    {\left(2-x^2\right)\left(4-x^2\right)\left(l^2-x^2\right)\sqrt{4l^2-4-l^2x^2} }
\text{\rm \ and \ }$$
$$\deriv{\tau} =\frac{\left(2-l^2\right)x }{\left(l^2-x^2\right)\sqrt{4l^2-4-l^2x^2 } }.$$
From these using $1\le x < l <\sqrt{2}$ we easily obtain that $\deriv{(\rho +\tau)} >0$ if and only if
$$x^4+\frac{10-7l^2 }{l^2-1 }x^2+\frac{10l^2-16 }{l^2-1 }<0.$$
Therefore $l=\sqrt{\frac{2d}{d+1}}$ yields that  $\deriv{(\rho +\tau)} >0$ if and only if
$$f(x)=x^4-\frac{4d-10}{d-1}x^2+\frac{4d-16}{d-1}<0.$$
Notice that $f(\sqrt{\frac{2(d-4)}{d-1} })=f(\sqrt{2})=0$ and $f(z)<0$ for all $\sqrt{\frac{2(d-4)}{d-1}}< z
< \sqrt{2}$. Finally, as $2\le \Vert {\bold x}_{d-1}-{\bold x}_0\Vert$ an easy computation shows that 
in fact, $x\le \sqrt{\frac{2(d-2)}{d-1}}$ and so 
$$\sqrt{\frac{2(d-4)}{d-1}} < \sqrt{ \frac{2(d-3)}{d-2}}\le x \le \sqrt{\frac{2(d-2)}{d-1} }< \sqrt{ \frac{2d}{d+1}} <\sqrt{2}.$$
Thus, $\rho +\tau$ as a function of $x$ is stricly increasing on the intervall 
$$\left[\sqrt{ \frac{2(d-3)}{d-2}}, \sqrt{\frac{2(d-2)}{d-1}} \right].$$ Hence, $\pi - ( \rho +\tau )$ as a function
of $x$ is streactly decreasing on the same interval.
Consequently, $\varphi$ as the largest possible value of  $\pi - ( \rho +\tau )$ is attained at $x=\sqrt{ \frac{2(d-3)}{d-2}}$. This yields $\cos \varphi = 
\frac{\sqrt{2}}{3}\frac{2d-1}{\sqrt{d(d-1)}}$ finishing the proof of Sublemma 4. \hskip05cm \qed
\bigskip

This completes the proof of Lemma 4.  $\qed$
\bigskip
As an immediate corollary of Lemma 4 we get the following statement.

\proclaim {Corollary 2} Let $\overline{B}\cap F_2$ be the $2-$dimensional base of the type I truncated wedge $\overline{W}_I$ (resp., type II truncated wedge $\overline{W}_{II}$) in the Voronoi polytope $P\subset {\bold E^d}$ of 
dimension $d\ge 8$. Then the number of line segments of positive length in
$\text{\rm relbd}(\overline{B}\cap F_2)$ is at most $4$.
\endproclaim

\heading 3. The lemma of comparison and \\
the integral representation of the surface density \\ 
            in (truncated) wedges of type I and II
\endheading

Recall that $B=\{ {\bold x}\in {\bold E^d}\vert \ \text{dist}({\bold o}, {\bold x})=\Vert{\bold x}\Vert\le 1\}$ and let $$S=\{ {\bold x}\in {\bold E^d}\vert \ \text{dist}({\bold o}, {\bold x})=\Vert{\bold x}\Vert = 1\}.$$
Then let $H\subset {\bold E^d}$ be a hyperplane disjoint from the interior of the unit ball $B$ and let
$Q\subset H$ be an arbitrary $(d-1)-$dimensional compact convex set. If $[{\bold o}, Q]$ denotes the convex
cone $\text{\rm conv}(\{{\bold o}\}\cup Q)$ with apex ${\bold o}$ and base $Q$, then {\it the (volume) density $\delta ([{\bold o}, Q], B)$} of the unit ball B in the cone $[{\bold o}, Q]$ is defined as
$$ \delta ([{\bold o}, Q], B)=\frac{\text{\rm Vol}_d ([{\bold o}, Q]\cap B) }{\text{\rm Vol}_d ( [{\bold o},
Q]) },$$
where $\text{\rm Vol}_d (\dots )$ refers to the corresponding $d-$dimensional Euclidean volume measure (\cite {4}). It is natural to introduce the following very similar notion.
\subhead Definition 5
\endsubhead
{\it The surface density $\widehat\delta ([{\bold o}, Q], S) $} of the unit sphere $S$ in the convex cone 
$[{\bold o}, Q]$ with apex ${\bold o}$ and base $Q$ is defined by
$$\widehat\delta ([{\bold o}, Q], S)= \frac{\text{\rm SVol}_{d-1}([{\bold o}, Q]\cap S)}{\text{\rm Vol}_{d-1}(Q) },$$
where $\text{\rm SVol}_{d-1}(\dots )$ refers to the corresponding $(d-1)-$dimensional spherical volume
measure.
\bigskip
If $h=\text{\rm dist}({\bold o}, H)$, then clearly $h\cdot\delta ([{\bold o}, Q], B)=\widehat\delta ([{\bold o}, Q], S)$. We will need
the following statement the first part of which is due to Rogers (\cite{23}) and the second part of which has been recently proved in \cite{4, Theorem 1}.

\proclaim {Lemma 5 (Lemma of Comparison)} Let $U=\text{\rm conv}\{{\bold o}, {\bold u}_1, \dots , {\bold u}_d\}$ be a
$d-$di\-men\-si\-o\-nal orthoscheme in ${\bold E^d}$ and let $V=\text{\rm conv}\{{\bold o}, {\bold v}_1, \dots ,
{\bold v}_d\}$ be a $d-$dimensional simplex of ${\bold E^d}$ such that $\Vert {\bold v}_i\Vert=\text{\rm dist}
({\bold o}, \text{\rm conv}\{{\bold v}_i, {\bold v}_{i+1}, \dots , {\bold v}_d\})$ for all $1\le i\le d-1$. If
$1\le \Vert {\bold u}_i\Vert \le \Vert {\bold v}_i\Vert $ holds for all $1\le i\le d$, then
$$\align 
\delta ( U, B) \ge \delta ( V, B) \text{\rm \ and \ } \tag1 \\
\widehat\delta (U, S) \ge \widehat\delta (V, S).  \tag2 
\endalign$$
\endproclaim 
 
At this point it is useful to introduce the following notations. (Notice that Sublemma 2 provides the
necessary geometry for Definition 7.)

\subhead Definition 6
\endsubhead
Let ${\bold x}_1, \dots , {\bold x}_n, n\ge 1$ be points in ${\bold E^d}, d\ge 1$ and let $X\subset {\bold E^d}$ be an arbitrary convex set. If $X_0=X$ and $X_m=\text{\rm conv}(\{{\bold x}_{n-(m-1)}\}\cup X_{m-1})$ for $m=1, \dots , n$, then we denote the final convex set $X_n$ by
$$[{\bold x}_1, \dots , {\bold x}_n, X].$$
\subhead Definition 7
\endsubhead
Let $W_I$ (resp., $\overline{W}_I$) denote the wedge (resp., truncated wedge) of type I
with the $2-$dimensional base $F_2$ (resp., $\overline{B}\cap F_2$) which is generated by the $(d-2)-$dimensional Rogers orthoscheme $\text{\rm conv}\{{\bold o}, {\bold r}_1, \dots , {\bold r}_{d-2}\}$ of the Voronoi polytope $P\subset {\bold E^d}, d\ge 4$. Then let
$$Q_I=[{\bold r}_1, \dots , {\bold r}_{d-3}, F_2] \text{\rm \ (resp., \ }\overline{Q}_I=[{\bold r}_1, \dots , {\bold r}_{d-3}, 
\overline{B}\cap F_2])$$
be called the {\it $(d-1)-$dimensional base} of the type $I$ wedge $W_I=[{\bold o}, Q_I]$ (resp., type $I$ truncated wedge $\overline{W}_I=
[{\bold o}, \overline{Q}_I]$). Similarly, we define the $(d-1)-$dimensional bases $Q_{II}$ and $\overline{Q}_{II}$ of $W_{II}$ and
$\overline{W}_{II}$. Finally, let 
$$h_1=\Vert{\bold r}_1\Vert, h_2=\Vert{\bold r}_2-{\bold r}_1\Vert, \dots , h_{d-2}=\Vert {\bold r}_{d-2}-{\bold r}_{d-3}\Vert.$$  

\bigskip
\proclaim { Sublemma 5}Let $W_I$ (resp., $W_{II}$) denote the wedge of type I (resp., of type $II$) 
with the $2-$dimensional base $F_2$ which is generated by the $(d-2)-$dimensional Rogers orthoscheme $\text{\rm conv}\{{\bold o}, {\bold r}_1, \dots , {\bold r}_{d-2}\}$ of the Voronoi polytope $P\subset {\bold E^d}, d\ge 4$. Then we have the following volume formulas.
\roster
\item $$\text{\rm Vol}_{d-1}(Q_I)=\frac{2}{(d-1)!}\left(\prod_{i=2}^{d-2}h_i\right)\text{\rm Vol}_2(F_2)\text{ \ and \ }$$
$$\text{\rm Vol}_d(W_I)=\frac{h_1}{d}\text{\rm Vol}_{d-1}(Q_I)=\frac{2}{d!}\left(\prod_{i=1}^{d-2}h_i\right)\text{\rm Vol}_2(F_2).$$
Similar formulas hold for the corresponding dimensional volumes of $\overline{Q}_I, \overline{W}_I,$ $Q_{II},$ $W_{II},$ $\overline{Q}_{II}$
and $\overline{W}_{II}$.
\item In general, if $K\subset\text{\rm aff}F_2$ is a convex domain, then
$$\text{\rm Vol}_{d-1}([{\bold r}_1, \dots , {\bold r}_{d-3}, K])=\frac{2}{(d-1)!}\left(\prod_{i=2}^{d-2}h_i\right)\text{\rm Vol}_2(K)\text{ \ and \ }$$
$$\text{\rm Vol}_d([{\bold o}, {\bold r}_1, \dots , {\bold r}_{d-3}, K )=\frac{h_1}{d}\text{\rm Vol}_{d-1}( [{\bold r}_1, \dots , {\bold r}_{d-3}, K])=\frac{2}{d!}\left(\prod_{i=1}^{d-2}h_i\right)\text{\rm Vol}_2(K).$$

\endroster
\endproclaim
\subhead Proof
\endsubhead
The proof follows from Sublemma 2 and Lemma 3 (Part (3)) in a straightforward way. \hskip0.5cm $\qed$
\bigskip
The central notion of this section is the limiting surface density introduced as follows (see also Fig. 7).

\OneFigure{file=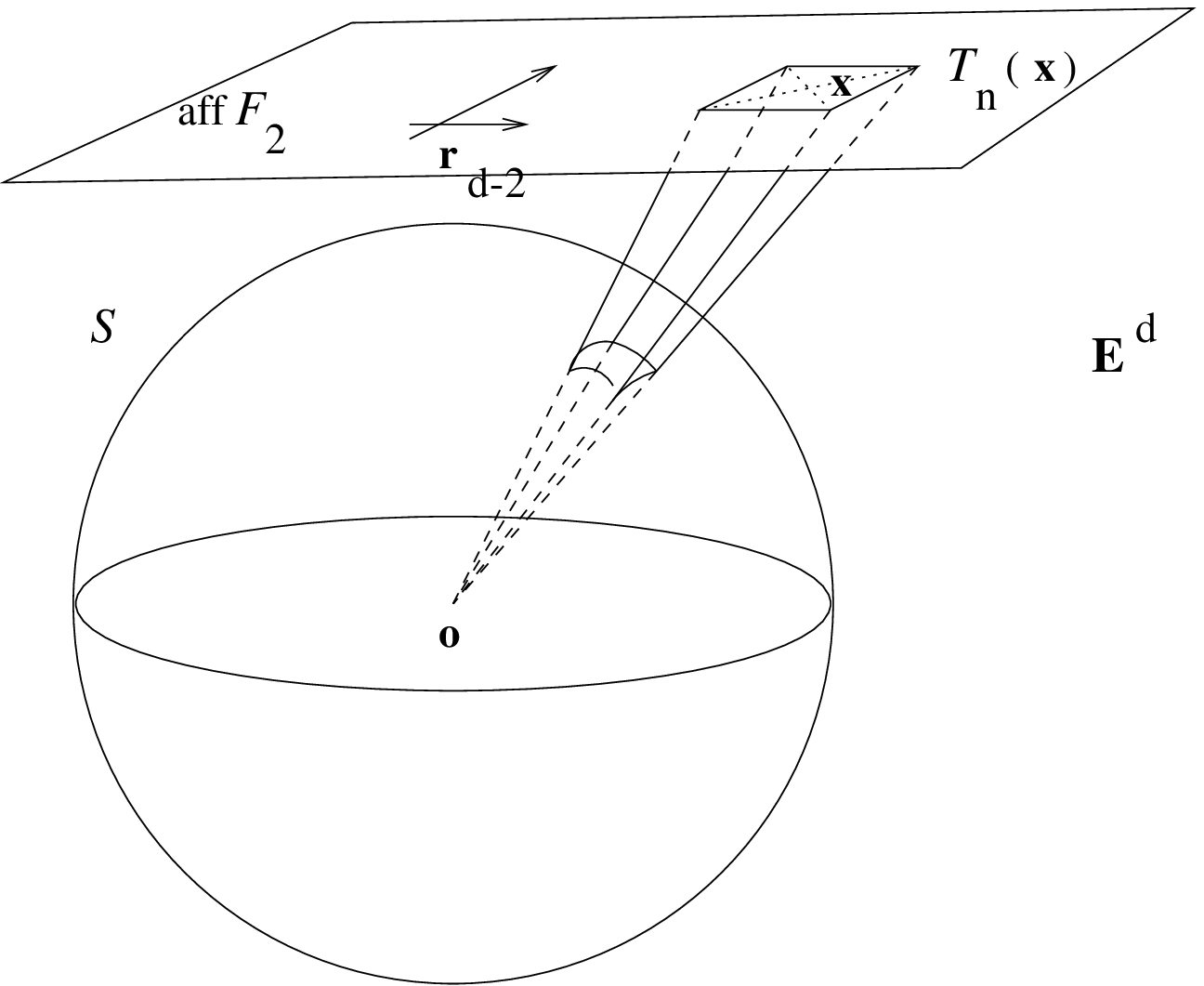}{Figure 7}

\subhead Definition 8
\endsubhead
Let $W_I$ (resp., $W_{II}$) denote the wedge of type I (resp., of type $II$) 
with the $2-$dimensional base $F_2$ which is generated by the $(d-2)-$dimensional Rogers orthoscheme $\text{\rm conv}\{{\bold o}, {\bold r}_1, \dots , {\bold r}_{d-2}\}$ of the Voronoi polytope $P\subset {\bold E^d}, d\ge 4$.
Then choose a coordinate system with two perpendicular axes in the plane $\text{\rm aff}F_2$ meeting at the
point ${\bold r}_{d-2}$. Now, if ${\bold x}$ is an arbitrary point of the plane $\text{\rm aff}F_2$, then for a positive integer $n$ let $T_n({\bold x})\subset\text{\rm aff}F_2$ denote the square centered at ${\bold x}$ having sides of 
length $\frac{1}{n}$ parallel to the fixed coordinate axes. Then {\it the limiting surface density} 
$ \widehat\delta_{\text{\rm lim}}\left([{\bold o}, {\bold r}_1, \dots , {\bold r}_{d-3}, {\bold x}], S\right) $
of the unit sphere $S$ in the $(d-2)-$dimensional orthoscheme $[{\bold o}, {\bold r}_1, \dots , {\bold r}_{d-3}, {\bold x}]$ is defined by
$$\widehat\delta_{\text{\rm lim}}\left([{\bold o}, {\bold r}_1, \dots , {\bold r}_{d-3}, {\bold x}], S\right)=
\lim_{n\to\infty}\widehat\delta\left([{\bold o}, {\bold r}_1, \dots , {\bold r}_{d-3}, T_n({\bold x})], S\right).$$
\bigskip


Based on this we are able to give an integral representation of the surface density in a (truncated) wedge.

\proclaim {Lemma 6}
Let $W_I$ (resp., $W_{II}$) denote the wedge of type I (resp., of type $II$) 
with the $2-$dimensional base $F_2$ which is generated by the $(d-2)-$dimensional Rogers orthoscheme $\text{\rm conv}\{{\bold o}, {\bold r}_1, \dots , {\bold r}_{d-2}\}$ of the Voronoi polytope $P\subset {\bold E^d}, d\ge 4$.
\roster
\item If ${\bold x}\in\text{\rm aff}F_2$ and ${\bold y}\in\text{\rm aff}F_2$ are points such that $\Vert {\bold x}\Vert \le \Vert {\bold y}\Vert$, then
$$\widehat\delta_{\text{\rm lim}}\left([{\bold o}, {\bold r}_1, \dots , {\bold r}_{d-3}, {\bold x}], S\right)\ge
\widehat\delta_{\text{\rm lim}}\left([{\bold o}, {\bold r}_1, \dots , {\bold r}_{d-3}, {\bold y}], S\right).$$
\item For the surface densities of the unit sphere $S$ in the wedge $W_I$ and in the truncated wedge $\overline{W}_I$ we have the following formulas:
$$\widehat\delta (W_I, S)=\frac{\text{\rm SVol}_{d-1}([{\bold o}, Q_I]\cap S)}{\text{\rm           Vol}_{d-1}(Q_I)}=$$
$$\frac{1}{ \text{\rm Vol}_2(F_2)}\int_{F_2}\widehat\delta_{\text{\rm lim}}\left([{\bold o}, {\bold r}_1, \dots , {\bold r}_{d-3}, {\bold x}], S\right)\ dx
\text{ \ and \ }$$
$$\widehat\delta (\overline{W}_I, S)=\frac{\text{\rm SVol}_{d-1}([{\bold o}, \overline{Q}_I]\cap S)}{\text{\rm Vol}_{d-1}(\overline{Q}_I)}=$$
$$\frac{1}{ \text{\rm Vol}_2(\overline{B}\cap F_2)}\int_{  \overline{B}\cap F_2}\widehat\delta_{\text{\rm lim}}\left([{\bold o}, {\bold r}_1, \dots , {\bold r}_{d-3}, {\bold x}], S\right)\ dx,$$
where $dx$ stands for the Euclidean area element in the plane $\text{\rm aff}F_2$. Similar formulas hold for $W_{II}$ and $\overline{W}_{II}$. 
\item In general, if $K\subset\text{\rm aff}F_2$ is a convex domain, then the surface density of the unit sphere $S$ in the $d-$dimensional convex cone $[{\bold o}, {\bold r}_1, \dots , {\bold r}_{d-3}, K]$
with apex ${\bold o}$ and $(d-1)-$dimensional base $[{\bold r}_1, \dots , {\bold r}_{d-3}, K]$
can be computed as follows:
$$\widehat\delta([{\bold o}, {\bold r}_1, \dots , {\bold r}_{d-3}, K], S)=
\frac{1}{ \text{\rm Vol}_2(K)}\int_{K}\widehat\delta_{\text{\rm lim}}\left([{\bold o}, {\bold r}_1, \dots , {\bold r}_{d-3}, {\bold x}], S\right)\ dx.$$

\endroster
\endproclaim 

\subhead Proof
\endsubhead  
\newline

(1) It is sufficient to look at the case $\Vert {\bold x}\Vert < \Vert {\bold y}\Vert$. (The case $\Vert {\bold x}\Vert = \Vert {\bold y}\Vert$ follows from this by standard limit procedure.) Then recall that $$\widehat\delta\left([{\bold o}, {\bold r}_1, \dots , {\bold r}_{d-3}, T_n({\bold x})], S\right)=h_1\delta \left([{\bold o}, {\bold r}_1, \dots , {\bold r}_{d-3}, T_n({\bold x})], S\right) 
\text{ \ and \ }$$
$$ \widehat\delta\left([{\bold o}, {\bold r}_1, \dots , {\bold r}_{d-3}, T_n({\bold y})], S\right)=h_1\delta \left([{\bold o}, {\bold r}_1, \dots , {\bold r}_{d-3}, T_n({\bold y})], S\right).$$
Thus, it is sufficient to show that if $n$ is sufficiently large, then
 
$$\delta \left([{\bold o}, {\bold r}_1, \dots , {\bold r}_{d-3}, T_n({\bold x})], S\right)\ge
\delta \left([{\bold o}, {\bold r}_1, \dots , {\bold r}_{d-3}, T_n({\bold y})], S\right).$$
This we can get as follows. We can approximate the $d-$dimensional convex cone 
$[{\bold o}, {\bold r}_1, \dots , {\bold r}_{d-3}, T_n({\bold x})]$ (resp., $[{\bold o}, {\bold r}_1, \dots , {\bold r}_{d-3}, T_n({\bold y})]$) arbitrarily close with a finite (but possible large) number of 
non-overlapping $d-$dimensional orthoschemes each containing the $(d-3)-$dimensional orthoscheme 
$[{\bold o}, {\bold r}_1, \dots , {\bold r}_{d-3}]$ as a face and each having all the edge lengths of
the $3$ edges 
going out from the vertex ${\bold o}$ and not lying on the face $[{\bold o}, {\bold r}_1, \dots , {\bold r}_{d-3}]$ close to $\Vert{\bold x}\Vert$ (resp., $\Vert {\bold y}\Vert$) for $n$
sufficiently large (see also Sublemma 2). Thus, the claim follows from Part (1) of Lemma 5 (Lemma of Comparison)
rather easily.   

(2)-(3) It is sufficient to prove the corresponding formula for $K$.

A typical term of the Riemann-Lebesgue sum of $$\frac{1}{ \text{\rm Vol}_2(K)}\int_{K}\widehat\delta_{\text{\rm lim}}\left([{\bold o}, {\bold r}_1, \dots , {\bold r}_{d-3}, {\bold x}], S\right)\ dx$$ is equal to
$$\frac{1}{ \text{\rm Vol}_2(K)}\widehat\delta\left([{\bold o}, {\bold r}_1, \dots , {\bold r}_{d-3}, T_n({\bold x}_m)], S\right)\text{\rm Vol}_2(T_n({\bold x}_m)), m\in M.$$
Using Sublemma 5 this turns out to be equal to
$$\frac{\text{\rm Vol}_{d-1}([{\bold r}_1, \dots , {\bold r}_{d-3}, T_n({\bold x}_m)]) } { \text{\rm Vol}_{d-1}([{\bold r}_1, \dots , {\bold r}_{d-3}, K])}\widehat\delta\left([{\bold o}, {\bold r}_1, \dots , {\bold r}_{d-3}, T_n({\bold x}_m)], S\right)=$$
$$\frac{\text{\rm SVol}_{d-1}([{\bold o}, {\bold r}_1, \dots , {\bold r}_{d-3}, T_n({\bold x}_m)]\cap S) }
{ \text{\rm Vol}_{d-1}([{\bold r}_1, \dots , {\bold r}_{d-3}, K]) }.$$
Finally, as the union of the nonoverlapping squares $T_n({\bold x}_m), m\in M$ is a good approximation
of the convex domain $K$ in the plane $\text{\rm aff}F_2$ we get that
$$\sum_{m\in M}^{ }\frac{\text{\rm SVol}_{d-1}([{\bold o}, {\bold r}_1, \dots , {\bold r}_{d-3}, T_n({\bold x}_m)]\cap S) }{ \text{\rm Vol}_{d-1}([{\bold r}_1, \dots , {\bold r}_{d-3}, K])}=$$
$$\frac{\sum_{m\in M}^{ } \text{\rm SVol}_{d-1}([{\bold o}, {\bold r}_1, \dots , {\bold r}_{d-3}, T_n({\bold x}_m)]\cap S)}{\text{\rm Vol}_{d-1}([{\bold r}_1, \dots , {\bold r}_{d-3}, K]) }$$
is a good approximation of 
$$\frac{\text{\rm SVol}_{d-1}\left([{\bold o}, {\bold r}_1, \dots , {\bold r}_{d-3}, K]\cap S\right) } { \text{\rm Vol}_{d-1}([{\bold r}_1, \dots , {\bold r}_{d-3}, K])}=\widehat{\delta}( [{\bold o}, {\bold r}_1, \dots , {\bold r}_{d-3}, K]  , S).$$ 
This completes the proof of Lemma 6. \hskip0.5cm $\qed$

\heading 4. Truncation of wedges increases the surface density
\endheading

\proclaim {Lemma 7}
Let $W_I$ (resp., $W_{II}$) denote the wedge of type I (resp., of type $II$) 
with the $2-$dimensional base $F_2$ which is generated by the $(d-2)-$dimensional Rogers orthoscheme $\text{\rm conv}\{{\bold o}, {\bold r}_1, \dots , {\bold r}_{d-2}\}$ of the Voronoi polytope $P\subset {\bold E^d}, d\ge 4$. Then
$$\widehat\delta (W_I, S)\le \widehat\delta (\overline{W}_I, S)\ \text{ \ (\ resp., \ }
\widehat\delta (W_{II}, S)\le \widehat\delta (\overline{W}_{II}, S)\text{\ )}.$$
\endproclaim
\subhead Proof
\endsubhead
Notice that Part (1) of Lemma 6 easily implies that if $0< \text{\rm Vol}_2(F_2\setminus\overline{B})$, then
for any ${\bold x}^*\in F_2$ with $\Vert {\bold x}^*\Vert=\sqrt{\frac{2d}{d+1}}$ we have that
$$\frac{1}{ \text{\rm Vol}_2(F_2\setminus\overline{B})}\int_{F_2\setminus\overline{B}}\widehat\delta_{\text{\rm lim}}\left([{\bold o}, {\bold r}_1, \dots , {\bold r}_{d-3}, {\bold x}], S\right)\ dx \ \le$$
$$\widehat\delta_{\text{\rm lim}}\left([{\bold o}, {\bold r}_1, \dots , {\bold r}_{d-3}, {\bold x}^*], S\right)\le$$
$$ \frac{1}{ \text{\rm Vol}_2(\overline{B}\cap F_2)}\int_{  \overline{B}\cap F_2}\widehat\delta_{\text{\rm lim}}\left([{\bold o}, {\bold r}_1, \dots , {\bold r}_{d-3}, {\bold x}], S\right)\ dx \ .$$
Thus, if $0< \text{\rm Vol}_2(F_2\setminus\overline{B})$, then Part (2) of Lemma 6 yields that 
$$\widehat\delta(W_I, S)=\frac{1}{ \text{\rm Vol}_2(F_2)}\int_{F_2}\widehat\delta_{\text{\rm lim}}\left([{\bold o}, {\bold r}_1, \dots , {\bold r}_{d-3}, {\bold x}], S\right)\ dx=$$
$$\frac{\text{\rm Vol}_2(\overline{B}\cap F_2) }{\text{\rm Vol}_2(F_2)}\cdot\frac{1}{\text{\rm Vol}_2(\overline{B}\cap F_2)}\int_{\overline{B}\cap F_2}\widehat\delta_{\text{\rm lim}}\left([{\bold o}, {\bold r}_1, \dots , {\bold r}_{d-3}, {\bold x}], S\right)\ dx \ +$$
$$\frac{\text{\rm Vol}_2(F_2\setminus\overline{B}) }{\text{\rm Vol}_2(F_2)}\cdot\frac{1}{\text{\rm Vol}_2(F_2\setminus\overline{B})}\int_{F_2\setminus\overline{B}}\widehat\delta_{\text{\rm lim}}\left([{\bold o}, {\bold r}_1, \dots , {\bold r}_{d-3}, {\bold x}], S\right)\ dx \ \le$$
$$\frac{\text{\rm Vol}_2(\overline{B}\cap F_2) }{\text{\rm Vol}_2(F_2)}\cdot\frac{1}{\text{\rm Vol}_2(\overline{B}\cap F_2)}\int_{\overline{B}\cap F_2}\widehat\delta_{\text{\rm lim}}\left([{\bold o}, {\bold r}_1, \dots , {\bold r}_{d-3}, {\bold x}], S\right)\ dx \ +$$
$$\frac{\text{\rm Vol}_2(F_2\setminus\overline{B}) }{\text{\rm Vol}_2(F_2)}\cdot \frac{1}{\text{\rm Vol}_2(\overline{B}\cap F_2)}\int_{\overline{B}\cap F_2} \widehat\delta_{\text{\rm lim}}\left([{\bold o}, {\bold r}_1, \dots , {\bold r}_{d-3}, {\bold x}], S\right)\ dx \ =$$
$$\frac{1}{\text{\rm Vol}_2(\overline{B}\cap F_2)}\int_{\overline{B}\cap F_2} \widehat\delta_{\text{\rm lim}}\left([{\bold o}, {\bold r}_1, \dots , {\bold r}_{d-3}, {\bold x}], S\right)\ dx \ =\widehat\delta(\overline{W}_I, S).$$
As the same method works for $W_{II}$ and $\overline{W}_{II}$ this completes the proof of Lemma 7.
\hskip0.5cm $\qed$

\heading 5. Maximum surface density in truncated wedges of type I and II
\endheading

\subhead The case of truncated wedges of type I
\endsubhead
Let $\overline{W}_I$ denote the truncated wedge of type I with the $2-$dimensional base $\overline{B}\cap F_2$ which is generated by the $(d-2)-$dimensional Rogers orthoscheme $\text{\rm conv}\{{\bold o}, {\bold r}_1, \dots , {\bold r}_{d-2}\}$ of the Voronoi polytope $P\subset {\bold E^d}, d\ge 8$. By assumption $F_2$ is a $2-$dimensional face of the Voronoi polytope $P$ with 
$$\sqrt{\frac{2(d-2)}{d-1} }\le h=R(F_2) < \sqrt{\frac{2(d-1)}{d} }.$$
Let $G_0\subset\text{\rm aff}F_2$ (resp., $G\subset\text{\rm aff}F_2$) denote the closed circular disc of
radius $g_0(h)=\sqrt{\frac{2d}{d+1}-h^2}$ (resp., $g(h)=\frac{2-h^2}{\sqrt{4-h^2}}$) centered at the
point ${\bold r}_{d-2}$. It is easy to see that $G\subset\text{\rm relint}G_0$ for all $ \sqrt{\frac{2(d-2)}{d-1} }\le h < \sqrt{\frac{2(d-1)}{d}}$. (Moreover, $G=G_0$ for 
$h=\sqrt{\frac{2(d-1)}{d}}$.) Notice that $G_0=\overline{B}\cap \text{\rm aff}F_2$ thus, Corollary 1 implies
that there is no vertex of the face $F_2$ belonging to the relative interior of $G_0$ (Fig. 8). Moreover, as
$h=R(F_2)<\sqrt{2}$ Lemma 1 yields that $\frac{2}{\sqrt{4-h^2}}\le R(F_1)$ holds for any side $F_1$ of
the face $F_2$ hence, $G\subset F_2$ and of course, $G\subset \overline{B}\cap F_2=G_0\cap F_2$. Now, let
$M\subset\text{\rm aff}F_2$ be a square circumscribed about $G$. A straightforward computation yields that $\frac{g_0(h)}{g(h)}$ is
a strictly decreasing function on the interval $\left[\sqrt{\frac{2(d-2)}{d-1} }  , \sqrt{\frac{2(d-1)}{d}}\right)$
(i.e. $\derive{\left(\frac{g_0(h)}{g(h)}\right)} <0$ on the interval $\left( \sqrt{\frac{2(d-2)}{d-1} } , \sqrt{\frac{2(d-1)}{d}} \right)$) and $\frac{g_0\left(\sqrt{\frac{2(d-2)}{d-1}}\right)}{g\left(\sqrt{\frac{2(d-2)}{d-1}}\right)}=\sqrt{\frac{2d}{d+1}}<\sqrt{2}$.
Thus, the vertices of the square $M$ do not belong to $G_0$. Finally, as $d\ge 8$
Corollary 2 implies that there are at most $4$ sides of the face $F_2$ that intersect the relative interior
of $G_0$.

\OneFigure{file=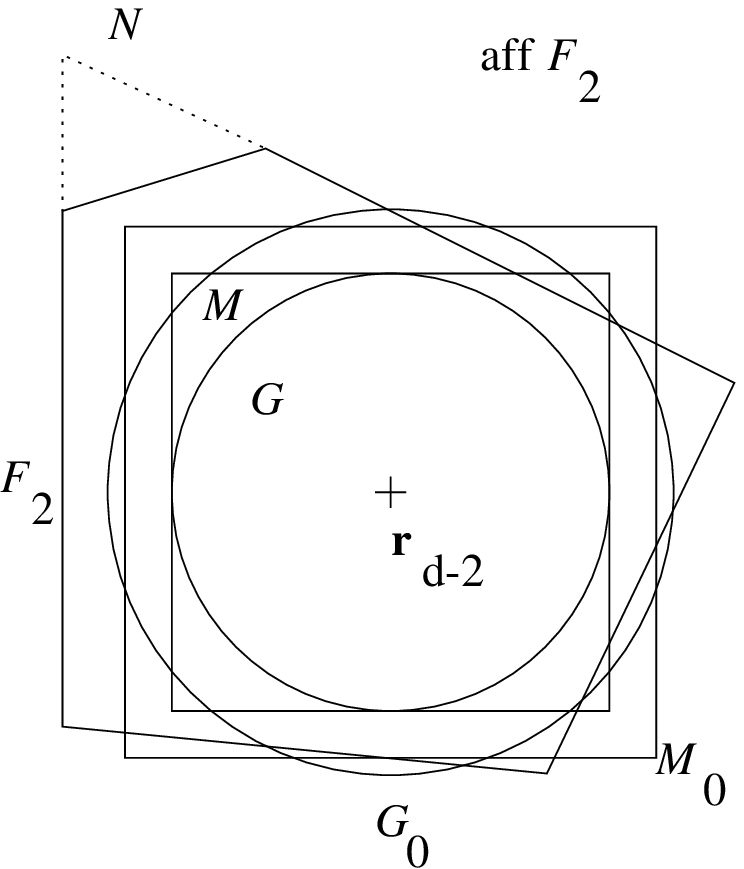}{Figure 8}

The following statement is the core part of this section.

\proclaim {Lemma 8} Let $\overline{W}_I$ denote the truncated wedge of type I with the $2-$dimensional base $\overline{B}\cap F_2$ which is generated by the $(d-2)-$dimensional Rogers orthoscheme $\text{\rm conv}\{{\bold o}, {\bold r}_1, \dots , {\bold r}_{d-2}\}$ of the Voronoi polytope $P\subset {\bold E^d}, d\ge 8$. Then
$$\widehat\delta(\overline{W}_I, S)\le \widehat\delta([{\bold o}, {\bold r}_1, \dots , {\bold r}_{d-3}, G_0\cap M], S).$$
\endproclaim

\subhead Proof
\endsubhead
Recall that according to Part (2) of Lemma 6 
$$\widehat\delta (\overline{W}_I, S)=\frac{1}{ \text{\rm Vol}_2(\overline{B}\cap F_2)}\int_{  \overline{B}\cap F_2}\widehat\delta_{\text{\rm lim}}\left([{\bold o}, {\bold r}_1, \dots , {\bold r}_{d-3}, {\bold x}], S\right)\ dx \ .$$
Moreover, Corollary 2 guarantees that the number of sides (i.e. of line segments of positive lengths)
of $\overline{B}\cap F_2=G_0\cap F_2$ is at most $4$. Thus, all these facts and Part (1) of Lemma 6 
imply that without loss of generality we may assume
that there exists a convex quadrangle $N\subset\text{\rm aff}F_2$ with 
$$G_0\cap F_2=G_0\cap N.$$
Now, if $G_0\cap N\ne G_0$ (resp., $ G_0\cap N= G_0$), then let $M_0\subset\text{\rm aff}F_2$ be a square (resp., a smallest square) centered at ${\bold r}_{d-2}$ with the property that
$$\text{\rm Vol}_2( G_0\cap N )=\text{\rm Vol}_2(G_0\cap M_0).$$
Obviously, no vertex of $M_0$ belongs to $G_0$.

\proclaim {Sublemma 6}
$$ \widehat\delta(\overline{W}_I, S)=\widehat\delta([{\bold o}, {\bold r}_1, \dots , {\bold r}_{d-3}, G_0\cap N], S)\le
\widehat\delta([{\bold o}, {\bold r}_1, \dots , {\bold r}_{d-3}, G_0\cap M_0], S).$$
\endproclaim
\subhead Proof
\endsubhead 
Take all possible convex  quadrilaterals $N'\subset\text{\rm aff}F_2$ with the property that no vertex of $N'$
belongs to the relative interior of $G_0$ and $G\subset N'$ moreover, $\text{\rm Vol}_2(G_0\cap N')=
\text{\rm Vol}_2(G_0\cap N)$. 

\OneFigure{file=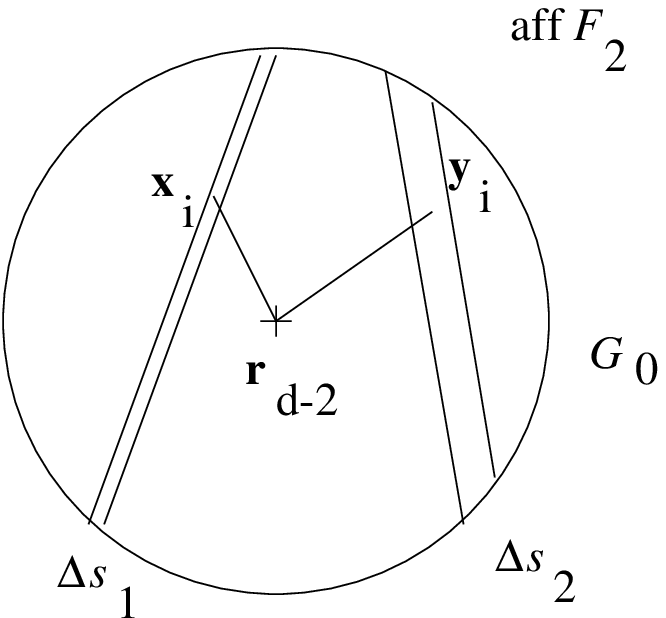}{Figure 9}

Obviously, there is convex quadrilateral of this family say, $N'$ for which
$$ \widehat\delta([{\bold o}, {\bold r}_1, \dots , {\bold r}_{d-3}, G_0\cap N'], S)=$$
$$\frac{1}{ \text{\rm Vol}_2(G_0\cap N')}\int_{G_0\cap N'}\widehat\delta_{\text{\rm lim}}\left([{\bold o}, {\bold r}_1, \dots , {\bold r}_{d-3}, {\bold x}], S\right)\ dx  $$
is maximal (see also Part (3) of Sublemma 6). We claim that all sides of $N'$ lie at the same distance from
${\bold r}_{d-2}$. We prove this by contradiction. Assume that there are two sides $s_1'$ and $s_2'$ of $G_0\cap N'$
such that the length of $s_1'$ is larger than the length of $s_2'$ (Fig. 9). Without loss of generality we
may assume that the endpoints of $s_1'$ and $s_2'$ are not vertices of $N'$. Then we move $s_1'$ 
(resp., $s_2'$) farther away from ${\bold r}_{d-2}$ (resp., closer to ${\bold r}_{d-2}$) by a small amount to get the new side $\widetilde{s}_1$ (resp., $\widetilde{s}_2$) such that the parallel strips $\varDelta s_1$ and $\varDelta s_2$ 
in $G_0$ determined by $s_1', \widetilde{s}_1$ and $s_2', \widetilde{s}_2$ have the same area. Then let $\widetilde{N}$ be the new convex quadrilateral obtained from $N'$ in the above manner. Now, it is easy to show that we can partition $\varDelta s_1$
(resp., $\varDelta s_2$) into an arbitrary large number say, $n$ of equal area convex subregions in each picking a point ${\bold x}_i, 1\le i\le n$ (resp., ${\bold y}_i, 1\le i\le n$) such that $\Vert {\bold x}_i-{\bold r}_{d-2}
\Vert < \Vert {\bold y}_i-{\bold r}_{d-2}\Vert$ i.e. $\Vert {\bold x}_i\Vert < \Vert {\bold y}_i\Vert$ for all $1\le i\le n$. As a result the proof of Part (1) of Lemma 6 yields the following inequality (see also the Lemma of Strict Comparison in Section 7):
$$ \widehat\delta_{\text{\rm lim}}\left([{\bold o}, {\bold r}_1, \dots , {\bold r}_{d-3}, {\bold x}_i], S\right) >
\widehat\delta_{\text{\rm lim}}\left([{\bold o}, {\bold r}_1, \dots , {\bold r}_{d-3}, {\bold y}_i], S\right)
\text{\rm \ for\  all\ } 1\le i\le n.$$
Hence, Part (3) of Lemma 6 implies in a straightforward way that
$$ \widehat\delta([{\bold o}, {\bold r}_1, \dots , {\bold r}_{d-3}, G_0\cap \widetilde{N}], S) >
\widehat\delta([{\bold o}, {\bold r}_1, \dots , {\bold r}_{d-3}, G_0\cap N'], S) \ ,$$
a contradiction. Thus, indeed all sides of $N'$ must lie at the same distance from ${\bold r}_{d-2}$ and as
a result we get via Lemma 6 that
$$ \widehat\delta([{\bold o}, {\bold r}_1, \dots , {\bold r}_{d-3}, G_0\cap N], S)\le
\widehat\delta([{\bold o}, {\bold r}_1, \dots , {\bold r}_{d-3}, G_0\cap M_0], S) ,$$
finishing the proof of Sublemma 6.  \hskip1.0cm $\qed$

\proclaim { Sublemma 7}
$$ \widehat\delta([{\bold o}, {\bold r}_1, \dots , {\bold r}_{d-3}, G_0\cap M_0], S)\le
 \widehat\delta([{\bold o}, {\bold r}_1, \dots , {\bold r}_{d-3}, G_0\cap M], S)  .$$
\endproclaim
\subhead Proof
\endsubhead 
Without loss of generality we may assume that $M$ and $M_0$ are homothetic with respect to ${\bold r}_{d-2}$ 
and $M\subset\text{\rm relint} M_0$. Let ${\bold u}\in G_0\cap M, {\bold v}\in G_0\cap M_0$
be midpoints of two corresponding and parallel sides of $G_0\cap M$ and $G_0\cap M_0$ (Fig. 10). Then let 
${\bold x}\in G_0\cap M$ and ${\bold y}\in G_0\cap M_0$ be the vertices of the above two corresponding and
parallel sides of $G_0\cap M$ and $G_0\cap M_0$ lying on the same side of the line ${\bold u}{\bold v}$ in
$\text{\rm aff}F_2$. Let ${\bold z}\in \text{\rm relbd}(G_0\cap M)$ be the point lying on the same side
of the line ${\bold u}{\bold v}$ in $\text{\rm aff}F_2$ as the points ${\bold x}, {\bold y}$ such that
$\angle{\bold u}{\bold r}_{d-2}{\bold z}={\bold v}{\bold r}_{d-2}{\bold z}=\frac{\pi}{4}$.
Let $U=\text{\rm conv}\{{\bold r}_{d-2}, {\bold u}, {\bold x}\}$, $V=\text{\rm conv}\{{\bold r}_{d-2}, {\bold v}, 
{\bold y}\}$. Moreover, let $X$ (resp., $Y$) be the circular sector of $G_0$ spanned by the center ${\bold r}_{d-2}$
and the shorter circular arc of $\text{\rm relbd}G_0$ between the points ${\bold x}, {\bold z}$ (resp., ${\bold y},
{\bold z}$). Now, Part (1) of Sublemma 2, Part (2) of Lemma 5 (Lemma of Comparison) and Lemma 6 easily yield
$$\widehat\delta([{\bold o}, {\bold r}_1, \dots , {\bold r}_{d-3}, {\bold r}_{d-2}, {\bold u} , {\bold x}], S)=
\widehat\delta([{\bold o}, {\bold r}_1, \dots , {\bold r}_{d-3}, U], S)\ge$$
$$\widehat\delta([{\bold o}, {\bold r}_1, \dots , {\bold r}_{d-3}, {\bold r}_{d-2}, {\bold v} , {\bold y}], S)=
\widehat\delta([{\bold o}, {\bold r}_1, \dots , {\bold r}_{d-3}, V], S)\ge$$
$$\widehat\delta([{\bold o}, {\bold r}_1, \dots , {\bold r}_{d-3}, X], S)=
\widehat\delta([{\bold o}, {\bold r}_1, \dots , {\bold r}_{d-3}, Y], S) .$$

\OneFigure{file=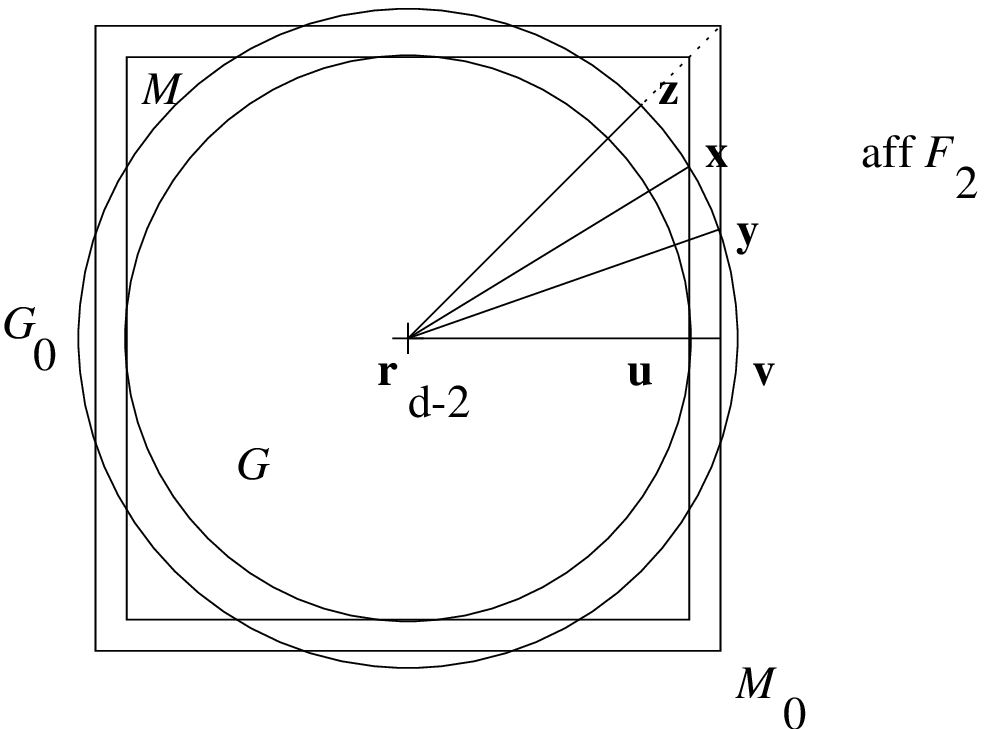}{Figure 10}

(Of course, in the case when $M_0$ is circumscribed about $G_0$ the points ${\bold y}, {\bold v}$ coincide and so,
in the above inequalities one has to replace $\widehat\delta([{\bold o}, {\bold r}_1, \dots , {\bold r}_{d-3}, V], S)$
by $\widehat\delta_{\text{\rm lim}}([{\bold o}, {\bold r}_1, \dots , {\bold r}_{d-3}, V], S) $.) Hence, these inequalities and the following obvious inequality (partly based on Sublemma 5)
$$\frac{\text{\rm Vol}_{d-1}([{\bold r}_1, \dots , {\bold r}_{d-3}, X])}{\text{\rm Vol}_{d-1}([{\bold r}_1, \dots , {\bold r}_{d-3}, U])}=\frac{\text{\rm Vol}_2(X)}{\text{\rm Vol}_2(U)}     \le
 \frac{\text{\rm Vol}_2(Y)}{\text{\rm Vol}_2(V)}= 
\frac{\text{\rm Vol}_{d-1}([{\bold r}_1, \dots , {\bold r}_{d-3}, Y])}{\text{\rm Vol}_{d-1}([{\bold r}_1, \dots , {\bold r}_{d-3}, V])}$$
imply (via Part (2) of Sublemma 2) in a straightforward way that
$$\widehat\delta([{\bold o}, {\bold r}_1, \dots , {\bold r}_{d-3}, U\cup X], S)=$$
$$\frac{ \text{\rm Vol}_{d-1}([{\bold r}_1, \dots , {\bold r}_{d-3}, U])}{ \text{\rm Vol}_{d-1}([{\bold r}_1, \dots , {\bold r}_{d-3}, U\cup X])} \widehat\delta([{\bold o}, {\bold r}_1, \dots , {\bold r}_{d-3}, U], S)+$$
$$\frac{ \text{\rm Vol}_{d-1}([{\bold r}_1, \dots , {\bold r}_{d-3}, X])}{ \text{\rm Vol}_{d-1}([{\bold r}_1, \dots , {\bold r}_{d-3}, U\cup X])} \widehat\delta([{\bold o}, {\bold r}_1, \dots , {\bold r}_{d-3}, X], S)\ge$$
$$\frac{ \text{\rm Vol}_{d-1}([{\bold r}_1, \dots , {\bold r}_{d-3}, V])}{ \text{\rm Vol}_{d-1}([{\bold r}_1, \dots , {\bold r}_{d-3}, V\cup Y])} \widehat\delta([{\bold o}, {\bold r}_1, \dots , {\bold r}_{d-3}, V], S)+$$
$$\frac{ \text{\rm Vol}_{d-1}([{\bold r}_1, \dots , {\bold r}_{d-3}, Y])}{ \text{\rm Vol}_{d-1}([{\bold r}_1, \dots , {\bold r}_{d-3}, V\cup Y])} \widehat\delta([{\bold o}, {\bold r}_1, \dots , {\bold r}_{d-3}, Y], S)=   $$
$$ \widehat\delta([{\bold o}, {\bold r}_1, \dots , {\bold r}_{d-3}, V\cup Y], S).  $$
Hence, by symmetry the inequality 
$$\widehat\delta([{\bold o}, {\bold r}_1, \dots , {\bold r}_{d-3}, G_0\cap M], S)\ge
 \widehat\delta([{\bold o}, {\bold r}_1, \dots , {\bold r}_{d-3}, G_0\cap M_0], S) $$
follows, finishing the proof of Sublemma 7. \hskip1.0cm $\qed$
\bigskip
Thus, Sublemma 6 and Sublemma 7 imply that
$$ \widehat\delta(\overline{W}_I, S)=\widehat\delta([{\bold o}, {\bold r}_1, \dots , {\bold r}_{d-3}, G_0\cap N], S)\le
\widehat\delta([{\bold o}, {\bold r}_1, \dots , {\bold r}_{d-3}, G_0\cap M_0], S)\le$$
$$\widehat\delta([{\bold o}, {\bold r}_1, \dots , {\bold r}_{d-3}, G_0\cap M], S).$$

This completes the proof of Lemma 8. \hskip1.0cm $\qed$
\bigskip
It is clear from the construction that we can write $ \widehat\delta([{\bold o}, {\bold r}_1, \dots , {\bold r}_{d-3}, G_0\cap M], S)$ as a function of $d-2$ variables namely,
$$\widehat\Delta (\xi_1, \dots , \xi_{d-3}, \xi_{d-2})=\widehat\delta([{\bold o}, {\bold r}_1, \dots , {\bold r}_{d-3}, G_0\cap M], S),$$
where $\xi_1=\Vert {\bold r}_1\Vert, \dots , \xi_{d-3}=\Vert{\bold r}_{d-3}\Vert, \xi_{d-2}=\Vert{\bold r}_{d-2}\Vert=h$. Corollary 1 and the assumption on $h$ imply that
$$m_1=1\le \xi_1, \dots , m_i=\sqrt{\frac{2i}{i+1}}\le \xi_i , \dots , m_{d-3}=\sqrt{\frac{2(d-3)}{d-2}}\le \xi_{d-3},$$
$$m_{d-2}=\sqrt{\frac{2(d-2)}{d-1} }\le \xi_{d-2}=h < \sqrt{\frac{2(d-1)}{d} }.$$
Notice that if $\Vert{\bold r}_i\Vert=m_i$ for all $1\le i\le d-2$, then $[{\bold o}, {\bold r}_1, \dots , {\bold r}_{d-3}, G_0\cap M]$ can be dissected into 4 pieces each being congruent to $W$ and therefore
$\widehat\delta([{\bold o}, {\bold r}_1, \dots , {\bold r}_{d-3},$ $G_0\cap M], S)=\widehat{\sigma}_d$.

\proclaim { Lemma 9}
$\widehat\Delta (\xi_1, \dots , \xi_{d-3}, \xi_{d-2})\le \widehat\Delta(m_1, \dots , m_{d-3}, m_{d-2})=
\widehat{\sigma}_d.$
\endproclaim
\subhead Proof
\endsubhead
For any fixed $\xi_{d-2}=h$ Part (2) of Lemma 5 (Lemma of Comparison) easily implies that
$$ \widehat\Delta (\xi_1, \dots , \xi_{d-3}, h)\le \widehat\Delta(m_1, \dots , m_{d-3}, h).$$
Finally, following the monotonicity idea of the proof of Sublemma 7 it is easy to show that the
function $\widehat\Delta(m_1, \dots , m_{d-3}, h)$ as a function of $h$ is decreasing on the interval
$[\sqrt{\frac{2(d-2)}{d-1} }, \sqrt{\frac{2(d-1)}{d} })$. (The proof is essentially based on the
fact that $\Vert {\bold u}\Vert, \Vert {\bold x}\Vert, \Vert{\bold z}\Vert, 
\frac{\text{\rm Vol}_2(X)}{\text{\rm Vol}_2(U)}$ are all increasing functions of $h$.) From this it
follows that
$$ \widehat\Delta(m_1, \dots , m_{d-3}, h)\le \widehat\Delta(m_1, \dots , m_{d-3}, m_{d-2})=
\widehat{\sigma}_d \ ,$$
finishing the proof of Lemma 9. \hskip1.0cm $\qed$
\bigskip
We conclude this section with the following immediate corollary of Lemma 8 and Lemma 9.

\proclaim {Corollary 3}
Let $\overline{W}_I$ denote the truncated wedge of type I with the $2-$dimensional base $\overline{B}\cap F_2$ which is generated by the $(d-2)-$dimensional Rogers orthoscheme $\text{\rm conv}\{{\bold o}, {\bold r}_1, \dots , {\bold r}_{d-2}\}$ of the Voronoi polytope $P\subset {\bold E^d}, d\ge 8$. Then
$$\widehat\delta(\overline{W}_I, S)\le \widehat{\sigma}_d.$$
\endproclaim

\subhead The case of truncated wedges of type II
\endsubhead
It is sufficient to prove the following statement.

\proclaim {Lemma 10}
Let $\overline{W}_{II}$ denote the truncated wedge of type II with the $2-$dimensional base $\overline{B}\cap F_2$ which is generated by the $(d-2)-$dimensional Rogers orthoscheme $\text{\rm conv}\{{\bold o}, {\bold r}_1, \dots , {\bold r}_{d-2}\}$ of the Voronoi polytope 
$P\subset {\bold E^d}, d\ge 4$. Then
$$ \widehat\delta(\overline{W}_{II}, S)\le \widehat{\sigma}_d.$$
\endproclaim
\subhead Proof
\endsubhead
By assumption $F_2$ is a $2-$dimensional face of the Voronoi polytope $P$ with 
$$\sqrt{\frac{2(d-1)}{d} }\le h=R(F_2) < \sqrt{\frac{2d}{d+1} }.$$
Let $G_0\subset\text{\rm aff}F_2$ denote the closed circular disc of radius $g_0(h)=\sqrt{\frac{2d}{d+1}-h^2}$ centered at the point ${\bold r}_{d-2}$. As $h=R(F_2)<\sqrt{2}$ Lemma 1
yields that
$$\sqrt{\frac{2d}{d+1}}\le \frac{2}{\sqrt{4-h^2}}\le R(F_1)$$
hold for any side $F_1$ of the face $F_2$. Thus, 
$$\overline{B}\cap F_2=G_0 \text{\rm \ and \ so\ }$$
$$ \widehat\delta(\overline{W}_{II}, S)=\widehat\delta([{\bold o}, {\bold r}_1, \dots , {\bold r}_{d-3}, 
G_0], S).$$
It is clear from the construction that we can write $ \widehat\delta([{\bold o}, {\bold r}_1, \dots , {\bold r}_{d-3}, G_0], S)$ as a function of $d-2$ variables namely,
$$\widehat{\Delta}^* (\xi_1, \dots , \xi_{d-3}, \xi_{d-2})=\widehat\delta([{\bold o}, {\bold r}_1, \dots , {\bold r}_{d-3}, G_0], S),$$
where $\xi_1=\Vert {\bold r}_1\Vert, \dots , \xi_{d-3}=\Vert{\bold r}_{d-3}\Vert, \xi_{d-2}=\Vert{\bold r}_{d-2}\Vert=h$. Corollary 1 and the assumption on $h$ imply that
$$m_1=1\le \xi_1, \dots , m_i=\sqrt{\frac{2i}{i+1}}\le \xi_i , \dots , m_{d-3}=\sqrt{\frac{2(d-3)}{d-2}}\le \xi_{d-3},$$
$$m_{d-2}^*=\sqrt{\frac{2(d-1)}{d} }\le \xi_{d-2}=h < \sqrt{\frac{2d}{d+1} }.$$
For any fixed $\xi_{d-2}=h$ Part (2) of Lemma 5 (Lemma of Comparison) easily implies that
$$ \widehat{\Delta}^* (\xi_1, \dots , \xi_{d-3}, h)\le \widehat{\Delta}^*(m_1, \dots , m_{d-3}, h).$$
Finally, applying again Part (2) of Lemma 5 we immediately get that
$$ \widehat{\Delta}^*(m_1, \dots , m_{d-3}, h)\le \widehat{\Delta}^*(m_1, \dots , m_{d-3}, m_{d-2}^*)\le
\widehat{\sigma}_d \ .$$
This completes the proof of Lemma 10. \hskip1.0cm $\qed$

\heading 6. Proof of the Theorem
\endheading
Let $P$ be a $d-$dimensional Voronoi polytope of a packing $\Cal P$ of $d-$dimensional unit balls
in ${\bold E^d}, d\ge 8$. Without loss of generality we may assume that the unit ball $B=\{ {\bold x}\in {\bold E^d}\vert \ \text{dist}({\bold o}, {\bold x})=\Vert{\bold x}\Vert\le 1\}$ centered at the origin ${\bold o}$ of $\bold E^d$ is one of the unit balls of $\Cal P$ with $P$ as its Voronoi cell. As before, let $S$ denote the boundary of $B$.

First, we dissect $P$ into $d-$dimensional Rogers simplices. Then let $\text{\rm conv}\{{\bold o}, {\bold r}_1,$ $\dots ,$ 
${\bold r}_d\}$ be one of these $d-$dimensional Rogers simplices assigned to the flag say, $F_0\subset \dots \subset F_{d-1}$ of $P$. As ${\bold r}_i\in F_{d-i}, 1\le i\le d$ it is clear that 
$\text{\rm aff}\{{\bold r}_{d-2}, {\bold r}_{d-1}, {\bold r}_d\}=\text{\rm aff}F_2$ and so
$$ \text{\rm dist}({\bold o}, \text{\rm aff}\{{\bold r}_{d-2}, {\bold r}_{d-1},
{\bold r}_d\})=\text{\rm dist}({\bold o}, \text{\rm aff}F_2)=R(F_2).$$
Notice that Corollary 1 implies that $\sqrt{\frac{2(d-2)}{d-1}}\le R(F_2)$. 

Second, we group the $d-$dimensional Rogers simplices of $P$ as follows.
\roster
\item If $ \sqrt{\frac{2(d-2)}{d-1}}\le R(F_2)<\sqrt{\frac{2(d-1)}{d}}$, then we assign the Rogers simplex $\text{\rm conv}\{{\bold o},$ ${\bold r}_1,$ $\dots , {\bold r}_d\}$ to the type I wedge $W_I$ with the $2-$dimensional base $F_2$
generated by the $(d-2)-$dimensional Rogers orthoscheme $\text{\rm conv}\{{\bold o}, {\bold r}_1, \dots , {\bold r}_{d-2}\}$ of the Voronoi polytope $P\subset{\bold E^d}, d\ge 8$.
\item If $\sqrt{\frac{2(d-1)}{d}}\le R(F_2) < \sqrt{\frac{2d}{d+1}}$, then we assign the Rogers simplex $\text{\rm conv}\{{\bold o},$ ${\bold r}_1,$ $\dots , {\bold r}_d\}$ to the type II wedge $W_{II}$ with the $2-$dimensional base $F_2$
generated by the $(d-2)-$dimensional Rogers orthoscheme $\text{\rm conv}\{{\bold o}, {\bold r}_1, \dots , {\bold r}_{d-2}\}$ of the Voronoi polytope $P\subset{\bold E^d}, d\ge 8$.
\item If $\sqrt{\frac{2d}{d+1}}\le R(F_2)$, then we assign the Rogers simplex  $\text{\rm conv}\{{\bold o}, {\bold r}_1, \dots , {\bold r}_d\}$ to itself as the type III wedge $W_{III}$.
\endroster
As the wedges of type I, II and III of the given Voronoi polytope $P$ sit over the $2-$skeleton of $P$ and form a tiling of $P$ it is clear that each $d-$dimensional Rogers simplex of $P$ belongs to exactly one of them. As a result, in order to show that the surface density  $\widehat\delta(P, S)=\frac{\text{\rm SVol}_{d-1}(S)}
{\text{\rm Vol}_{d-1}(\text{\rm bd}P) }=\frac{d\omega_d}{\text{\rm Vol}_{d-1}(\text{\rm bd}P) }$ of the unit sphere $S$ in the Voronoi polytope $P$ is bounded from above by $ \widehat{\sigma}_d$, it is sufficient to prove the following
inequalities 

\centerline{ $(\widehat{1})\hskip1.0cm \widehat\delta(W_I, S)\le  \widehat{\sigma}_d \ ;$}

\centerline{ $(\widehat{2})\hskip1.0cm \widehat\delta(W_{II}, S)\le  \widehat{\sigma}_d \ ;$}

\centerline{ $(\widehat{3})\hskip1.0cm \widehat\delta(W_{III}, S)\le  \widehat{\sigma}_d \ .$}

This final task left is now easy. Namely, Lemma 7, Corollary 3 and Lemma 10 yield $ (\widehat{1})$ and $(\widehat{2}) $
in a straightforward way. Finally, $(\widehat{3})$ follows with the help of Part (2) of Lemma 5 rather easily.

This completes the proof of the Theorem.

\heading 7. Proof of the Proposition 
\endheading
The proof of Lemma 5 published in \cite{4} (Part(2)) and in \cite{23} (Part(1)) can be modified in a
straightforward way such that it leads to the following somewhat stronger version of the Lemma of Comparison.

\proclaim {Lemma 11 (Lemma of Strict Comparison)} Let $U=\text{\rm conv}\{{\bold o}, {\bold u}_1, \dots , {\bold u}_d\}$ be a $d-$di\-men\-si\-o\-nal orthoscheme in ${\bold E^d}$ and let $V=\text{\rm conv}\{{\bold o}, {\bold v}_1, \dots ,{\bold v}_d\}$ be a $d-$dimensional simplex of ${\bold E^d}$ such that $\Vert {\bold v}_i\Vert=\text{\rm dist}
({\bold o}, \text{\rm conv}\{{\bold v}_i, {\bold v}_{i+1}, \dots , {\bold v}_d\})$ for all $1\le i\le d-1$. If
$1\le \Vert {\bold u}_i\Vert \le \Vert {\bold v}_i\Vert $ holds for all $1\le i\le d$ and there is an $i_0, 1\le i_0
\le d$ such that $1\le\Vert{\bold u}_{i_0}\Vert <\Vert{\bold v}_{i_0}\Vert$, then
$$\align 
\delta ( U, B) > \delta ( V, B) \text{\rm \ and \ } \tag1 \\
\widehat\delta (U, S) > \widehat\delta (V, S).  \tag2 
\endalign$$
\endproclaim 
First, recall that 
$$\widehat{\sigma}_d=\frac{\text{\rm Vol}_{d}(W\cap B)}{\text{\rm Vol}_{d}(W)}.$$
Second, an easy application of Part (1) of Lemma 11 implies that
$$ \lambda_d=\frac{\text{\rm Vol}_{d}\left( [{\bold o}, {\bold w}_1, \dots , {\bold w}_{d-3},
\vartriangleleft {\bold w}_{d-2}{\bold w}_{d}{\bold w}_{d+1}]\cap B\right)}
{\text{\rm Vol}_{d}\left([{\bold o}, {\bold w}_1, \dots , {\bold w}_{d-3},
\vartriangleleft {\bold w}_{d-2}{\bold w}_{d}{\bold w}_{d+1}]\right)} < $$ 
$$\frac{\text{\rm Vol}_{d}([{\bold o}, {\bold w}_1, \dots , {\bold w}_d]  \cap B)}{\text{\rm Vol}_{d}([{\bold o}, {\bold w}_1, \dots , {\bold w}_d] )}=\sigma_d.$$
Thus,
$$\widehat{\sigma}_d=\frac{\text{\rm Vol}_{d}(W\cap B)}{\text{\rm Vol}_{d}(W)}=$$
$$\frac{ \lambda_d  \text{\rm Vol}_{d}\left([{\bold o}, {\bold w}_1, \dots , {\bold w}_{d-3},
\vartriangleleft {\bold w}_{d-2}{\bold w}_{d}{\bold w}_{d+1}]\right)+\sigma_d  \text{\rm Vol}_{d}([{\bold o}, {\bold w}_1, \dots , {\bold w}_d] )} { \text{\rm Vol}_{d}(W) }<\sigma_d .$$

This completes the proof of the Proposition.

\heading Acknowledgement 
\endheading

I would like to thank Cornell University and especially Ro\-bert Connelly
for their hospitality while this work was undertaken. Also, the paper benefited from 
the valuable remarks of the three anonymous referees.


\bigskip


\Refs

\ref
\no 1
\by K. M. Ball
\pages 217-221
\paper A lower bound for the optimal density of lattice packings
\yr 1992
\vol 68
\jour Duke J. Mathematics
\endref

\ref
\no 2
\by E. Baranovskii
\pages 14-24
\paper On packing n-dimensional Euclidean space by equal spheres
\yr 1964
\issue 2
\vol 39
\jour Iz. Vissih Uceb. Zav. Mat.
\endref

\ref
\no 3
\by A. Bezdek, K. Bezdek and R. Connelly
\pages 111-130
\paper Finite and uniform stability of sphere packings
\yr 1998
\vol 20
\jour Discrete Comput. Geom.
\endref


\ref 
\no 4
\by K. Bezdek
\pages 131-143
\paper On a stronger form of Rogers's lemma and the minimum surface area
of Voronoi cells in unit ball packings
\yr 2000
\vol 518
\jour J. reine angew. Math.
\endref

\ref
\no 5
\by K. Bezdek and R. Connelly
\pages 185-200
\paper Two-distance preserving functions from Euclidean space
\yr 1999
\vol 39 
\issue 1-3
\jour Periodica Math. Hung.
\endref

\ref
\no 6
\by K. B\"or\"oczky
\pages 243-261
\paper Packing equal spheres in spaces of constant curvature
\yr 1978
\vol 32
\issue 3-4
\jour Acta Math. Acad. Sci. Hung.
\endref

\ref
\no 7
\by H. Cohn and N. Elkies
\pages 1-21
\paper New bounds on sphere packings I
\yr September 30, 2001
\jour Published electronically at http://ar Xiv.org/abs/math.MG/011009
\endref

\ref
\no 8
\by J. H. Conway and N. J. A. Sloane
\pages 593-620
\paper Laminated lattices
\yr 1982
\vol 116
\jour Ann. Math.
\endref

\ref
\no 9
\by J. H. Conway and N. J. A. Sloane
\pages 309-313
\paper The antipode construction for sphere packings
\yr 1996
\vol 123
\jour Invent. math.
\endref

\ref
\no 10
\by J. H. Conway and N. J. A. Sloane
\book Sphere Packings, Lattices and Groups
\publ Springer-Verlag, NY, 3rd edition
\yr 1998
\endref

\ref
\no 11
\by M. H. Dauenhauer and H. J. Zassenhaus
\pages 129-146
\paper Local optimality of critical lattice sphere-packing of regular tetrahedra
\yr 1987
\vol 64
\jour Discrete Math.
\endref


\ref
\no 12
\by G. Fejes T\'oth and W. Kuperberg
\book Packing and covering with convex sets, {\rm Handbook of Convex Geometry}
\eds P. M. Gruber and J. M. Wills
\publ Elsevier Science Publishers
\yr 1993
\page 799-860
\endref

\ref
\no 13
\by T. C. Hales
\pages 181-197
\paper Remarks on the density of sphere packings in three dimensions
\yr 1993
\vol 13
\jour Combinatorica
\endref


\ref
\no 14
\by T. C. Hales 
\pages 1-51
\paper Sphere packings I
\yr 1997
\vol 17
\jour Discrete Comput. Geom.
\endref

\ref
\no 15
\by T. C. Hales 
\pages 135-149
\paper Sphere packings II
\yr 1997
\vol 18
\jour Discrete Comput. Geom.
\endref

\ref
\no 16
\by T. C. Hales
\book The Kepler Conjecture
\publ Manuscript
\yr 1998
\endref

\ref
\no 17
\by W.-Y. Hsiang
\paper The optimal density of sphere packings in dimension $8$ and the uniqueness theorem on final packings with optimal density
\yr March 15, 2001
\jour Berkeley Mathematics Department Colloquium
\endref



\ref
\no 18
\by G. A. Kabatjanskii and V. I. Leven\v stein
\pages 3-25; English translation in Problems of Information Transmission ${\bold 14}$ (1978), 1-17
\paper Bounds for packings on a sphere and in space (in Russian)
\yr 1978
\vol 14
\jour Problemy Peredachy Informatsii
\endref

\ref
\no 19
\by F. R. Kschischang and S. Pasupathy
\pages 227-246
\paper Some ternary and quaternary codes and associated sphere packings
\yr 1992
\vol 38
\jour IEEE Trans. Inform. Theory
\endref

\ref
\no 20
\by J. Leech and N. J. A. Sloane
\pages 718-745
\paper Sphere packing and error-correcting codes
\yr 1971
\vol 23
\jour Canad. J. Math.
\endref



\ref
\no 21
\by D. J. Muder
\pages 351-375
\paper A new bound on the local density of sphere packings
\yr 1993
\vol 10
\jour Discrete Comput. Geom.
\endref

\ref
\no 22
\by G. Nebe and N. J. A. Sloane
\paper Table of densest packings presently known
\jour Published electronically at http://www.research.att.com/njas/lattices/density.html
\endref 


\ref
\no 23
\by C. A. Rogers
\pages 609-620
\paper The packing of equal spheres
\yr 1958
\vol 3
\issue 8
\jour Proc. London Math. Soc.
\endref

\ref
\no 24
\by C. A. Rogers
\book Packing and Covering
\publ Camb. Univ. Press
\yr 1964
\endref

\ref
\no 25
\by N. J. A. Sloane
\pages 387-396
\paper The sphere packing problem
\yr 1998
\vol III
\jour Doc. Math. J. DMV
\endref

\ref
\no 26
\by A. Vardy
\pages 119-133
\paper A new sphere packing in 20 dimensions
\yr 1995
\vol 121
\jour Invent. math.
\endref



\endRefs
\enddocument